# Formation control Part II: Algebraic aspects of information flow and singularities

M.-A. Belabbas [*]

January 14, 2011

## Abstract

Given an ensemble of autonomous agents and a task to achieve cooperatively, how much do the agents need to know about the state of the ensemble and about the task in order to achieve it? We introduce new methods to understand these aspects of decentralized control. Precisely, we introduce a framework to capture what agents with partial information can achieve by cooperating and illustrate its use by deriving results about global stabilization of directed formations. This framework underscores the need to differentiate the knowledge an agent has about the task to accomplish from the knowledge an agent has about the current state of the system.

The control of directed formations has proven to be more difficult than initially thought, as is exemplified by the lack of global result for formations with $n \geq 4$ agents. We established in part I that the space of planar formations has a non-trivial global topology. We propose here an extension of the notion of global stability which, because it acknowledges this non-trivial topology, can be applied to the study of formation control. We then develop a framework that reduces the question of whether feedback with partial information can stabilize the system to whether two sets of functions intersect. We apply this framework to the study of a directed formation with $n = 4$ agents and show that the agents do not have enough information to implement locally stabilizing feedback laws. Additionally, we show that feedback laws that respect the information flow cannot stabilize a target configuration without stabilizing other, unwanted configurations.

## Contents

---
[*]M.-A. Belabbas is with the School of Engineering and Applied Sciences, Harvard University, Cambridge, MA 02138 `belabbas@seas.harvard.edu`

# 1 Introduction

In Part I of this paper, we studied the geometry of the space of formations in the plane. In particular, we established that the space of edge lengths-normalized formations of $n$ agents is a complex projective space of dimension $n - 2$ and we described it explicitly for $n = 3$ and $n = 4$. This description relied on exhibiting discrete symmetry groups acting on frameworks and we conjectured the existence of a relation between these groups and Henneberg sequences. By relating the convexity of the space of edge lengths to a topological characteristic of $\mathbb{CP}(n - 2)$, we also provided a lower bound on the cardinality of the symmetry groups. In the last part, we presented a general dynamical model for formation control that respected both the invariance of the system under rigid transformations of the plane and the information flow as described by the underlying graph. Let us call this class of models $\mathcal{F}$.

In this part, we define two distinct global stability properties of formation control and investigate them in the case of the 2-cycles formation. Namely, we know from Part I that there are 4 frameworks in the plane that correspond to a generic vector of edge lengths for the 2-cycles. We answer the two following questions: "does there exist a feedback system in $\mathcal{F}$ that makes these four frameworks locally stable?" and "does there exist a control law that makes any of these four frameworks stable, and no other frameworks stable?" The first question thus asks for the existence of a globally defined control law that will make the four formations corresponding to a given edge lengths locally stable, whereas the second question is concerned with the existence of other stable formations.

We answer these questions by introducing new approaches to analyze the space $\mathcal{F}$. We first define an algebraic framework *to capture the range of behaviors that are realizable by elements of $\mathcal{F}$*. Loosely speaking, the main idea behind our framework is the following: we associate to elements of $\mathcal{F}$ a subset in the ring of smooth real-valued functions on the state-space, and to the objective that we want the decentralized control problem to achieve an ideal in that ring. We then reduce the feasibility question to an intersection problem between these two spaces. Within this framework, questions such as how much more information should an agent have in order to achieve an objective take a particularly natural form.

Indeed, we essentially reduce the problem of knowing whether an ensemble of agents can accomplish a certain task to problems of the following type: assume the task is encoded in the zero set of a function $p(x+y)$, and that the range of behaviors of the decentralized system for this task is encoded as the product $u_1(x)u_2(y)$. The feasibility of the objective is thus reduced to finding $u_1(x)$ and $u_2(y)$ smooth and such that

$$u_1(x)u_2(y) = 0 \Leftrightarrow p(x+y) = 0,$$

where it is understood here that agent 1 observes $x$ and agent 2 observes $y$. If, for example, $p$ is the identity, then this is clearly possible if and only if $u_1$ or $u_2$ are identically zero. Hence, we conclude that we would need to let agent 1 or 2 access more information (e.g. let $u_1$ depend on $y$) in order to solve this problem in a non-trivial manner. We revisit this example in Section 5, where these ideas are given a formal framework. An additional advantage of this framework is that it allows the use of the many computational tools that have been developed, e.g. Gröbner bases, to handle the algebraic structures involved.

The other method developed is the use of singularity theory to prove the existence of ancillary equilibria in the dynamics of decentralized control systems. With the exception of some work on the relation between Lyapunov theory and Morse theory [WY73], ideas from topology and singularity theory have not played an important role in the analysis of globally stabilizing control laws. In this paper, we will show that such ideas can prove to be fundamental in understanding how well one can hope to do in a control design task. In particular, we will use ideas from singularity theory [AAIS94, GSS88] to prove that if one tries to stabilize a system at a given configuration via a continuous feedback law, other equilibrium configurations appear and are stable, thus preventing global stability.

The 2-cycles formation (Figure 1b) is studied in depth in this paper. Its importance stems from its status as the second simplest problem in the class of directed formation problems, the simplest one being the triangular formation. This is a consequence of a theorem of Baillieul and Suri [BS03] that built upon earlier work of Brockett [Bro83]. The result asserts that when the interactions are asymmetric, or directed, one can generally require of an agent to have two leaders at the most. From this point of view, we can understand the formation depicted in Figure 1b as being the formation right above the triangle in terms of complexity. Indeed, this formation was singled out in [CAM$^+$10] as the prototypical example of the difficulty to make progress in this area and it underscored the need of new results to address decentralized control problems.

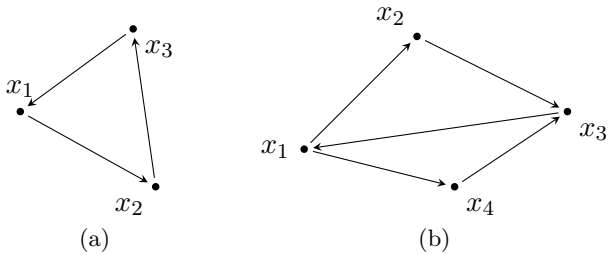

Figure 1: In Figure 1a, we show three agents in a cyclic formation in the plane. Agent 1 observes agent 2, which observes agent 3, which in turn observes agent 1. In Figure 1b, we depict the 2-cycles problem analyzed in detail in this paper. Agent 1 observes agents 2 and 4, agents 2 and 4 observe agent 3 which observes agent 1.

The results that are the closest to the formation control aspect of this work revolve around the triangle formation, local stabilization properties and formation control with undirected information flow. The triangular formation was dealt with in a series of papers ([CAM$^+$10, CMY$^+$07] and references therein). More precisely, it was shown that there exists a control law $u$ which will satisfactorily stabilize triangular formations. Furthermore, a careful analysis exhibited a whole family of such laws and gave a unified framework for the analysis of their convergence properties. In [YADF09], the local stabilization of a given minimally persistent [HADB07] formation was related to some easily checkable condition involving minors of a Jacobian. While the first result is global, and carefully addresses the existence of a set of initial conditions that may result in the system evolving to an undesired configuration, it does not address formations with more than 3 agents. In contrast, the latter result handles formations of arbitrary size, but only addresses local properties. For similar work on the control of formations with symmetric interactions, we refer the reader to [KBF08, SBF06].

In order to put our results in perspective, observe that the work of Yu, Anderson, Dasgupta and Fidan [YADF09] asserts that one can find a locally stabilizing control law, via choosing appropriate feedback gains for the linearized system, for any given configuration. The control laws were computed with complete knowledge of the formation, thus separating the design stage from the dynamics stage. Our stability results state that, at least in the case of the two-cycles, these local control laws *cannot be evaluated by the agents with the information they have at their disposal.* Moreover, we prove that such control laws introduce other stable equilibria.

The paper is organized as follows: in the second section, we introduce a notion of global stability that is more appropriate to the study of decentralized control problems. This notion of stability is weaker than the usual one and acknowledges some simple facts from Morse theory. In Section 3, we review the necessary background in singularity theory, relegating the most technical aspects of transversality and Thom's theorem to an appendix. In Section 4, we briefly present the logistic equation. We introduce in Section 5 an algebraic framework to capture the range of behaviors of decentralized control systems, and we apply it in the following section by proving that one cannot make all frameworks corresponding to a given edge length vector stable. In Section 7, we prove that feedback laws on the 2-cycles will generically make a stable ancillary equilibrium appear. In the last section, we discuss some forthcoming results and illustrate how a control law that stabilizes the triangle formation fails to stabilize the 2-cycles.

We state the results that relate to the 2-cycles explicitly here: let $x_1, x_2, x_3, x_4 \in \mathbb{R}^2$ represent the positions of the agents in the plane, and define

$$z_1 = x_2 - x_1, z_2 = x_3 - x_2, z_3 = x_1 - x_3, z_4 = x_3 - x_4, z_5 = x_3 - x_1. \tag{1}$$

Let
$$e_i = z_i \cdot z_i - d_i,$$
where
$$x \cdot y$$
is the inner product of vectors $x$ and $y$, be the error between the edge lengths and the target edge length $d_i$. The objective of the formation control problem is to find a control law that drives the system to a configuration with $e_i = 0$ for $i = 1, \ldots, 5$.

Consider the class of distributed feedback control systems, which respect the invariance under rigid transformations:
$$\begin{cases} \dot{x}_1 &= u_1 z_1 + u_5 z_5 \\ \dot{x}_2 &= u_2 z_2 \\ \dot{x}_3 &= u_3 z_3 \\ \dot{x}_4 &= u_4 z_4 \end{cases} \tag{2}$$

where $u_1$ and $u_5$ are smooth functions of $d_1, d_5, z_1, z_5$ and $u_i$ are smooth functions of $d_i$ and $z_i$ for $i = 2, 3, 4$. We will show that, for a set of positive measure of edge lengths, there are no *robust* distributed control laws that stabilize equilibria with edge lengths $(d_1, ..., d_5)$, and moreover, there will be a stable equilibrium with $e_i \neq 0$. In other words, if there exists a $u$ that accomplishes any of the two tasks above (i.e. stabilizing the four equilibria with $e_i = 0$, or stabilizing at least one such equilibria while making sure no other formation is stable), then there is an $r > 0$ such that for all $0 < \varepsilon < r$, and almost all $\tilde{u}_i$, the system

$$\begin{cases} \dot{x}_1 &= (u_1 + \varepsilon \tilde{u}_1)z_1 + (u_5 + \varepsilon \tilde{u}_5)z_5 \\ \dot{x}_2 &= (u_2 + \varepsilon \tilde{u}_2)z_2 \\ \dot{x}_3 &= (u_3 + \varepsilon \tilde{u}_3)z_3 \\ \dot{x}_4 &= (u_4 + \varepsilon \tilde{u}_4)z_4 \end{cases} \quad (3)$$

will fail to do both. Hence, the control law is not robust to the least error in modelling or measurement. An important aspect is that we need to only consider the set of perturbations that *respect the distributed nature* of the control law.

## 2 Notions of stability

Consider the control system
$$\dot{x} = f(x, u(x)) \quad (4)$$
where $x \in M$, a smooth manifold, and all functions are assumed smooth.

We are mostly interested in *global* results about stabilization. From Part I, we know that formation control problems evolve on a manifold $M$ with non-trivial homology groups [War83]. As a result of the Morse inequalities [Sma67], these systems cannot be globally stable in the usual sense: there is no continuous $u$ such that (4) has a *unique* equilibrium. Such situations happen frequently in nonlinear control, e.g. in steering control.

**Example 1.** *Control problems where one needs to stabilize the direction of motion of a vehicle are modelled by a dynamical system evolving on the unit circle $S^1$:*

$$\dot{x} = f(x), x \in S^1.$$

*In this one-dimensional case, the Morse inequalities reduce to the Poincaré-Hopf identity:*

$$\sum_{x=equilibria} \mathrm{ind}_\mathrm{x} f(x) = \chi(S^1),$$

*where $\chi(M)$ is the Euler characteristic of $M$ and $\mathrm{ind}_x$ is the index of the vector field at the equilibrium $x$. In the case of the circle, we can prove that exponentially stable equilibria correspond to equilibria of index $-1$ and, reciprocally, exponentially unstable equilibria to equilibria of index $1$. Since the Euler characteristic of $S^1$ is $\chi(S^1) = 0$, the Poincaré-Hopf identity becomes:*

$$\# \text{ exponentially stable equil.} = \# \text{ exponentially unstable equilibria }.$$

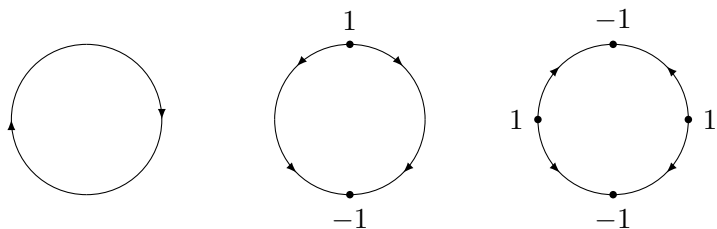

Figure 2: We represent three continuous vector fields on the circle. The one on the left has no equilibrium, the one in the center has one stable equilibrium and one unstable equilibrium and the one on the right has two stable equilibria and two unstable equilibria; the indices of the equilibria are indicated on the Figure. Because the sum of the indices is constrained to be zero, there is an even number of equilibria and, in particular, no continuous system on the circle can be globally stable.

*We illustrate it in Figure 2.*

*For every stable equilibrium, the topology of $S^1$ thus forces the appearance of an unstable equilibrium. Hence, system evolving on $S^1$ cannot be globally stable in the usual sense.*

From a practical point of view, however, if one could make one equilibrium stable, and all other equilibria either saddles or unstable, the system would behave as if it were globally stable. Indeed, a vanishingly small perturbation would ensure that the system, if at a saddle or unstable equilibrium, evolves to the unique stable equilibrium. We formalize and elaborate on this observation here.

Let $\mathcal{E}_d$ be a finite subset of $M$ containing configurations that we would like to stabilize via feedback. All configurations in $\mathcal{E}_d$ are equally appropriate for the stabilization purpose. We are thus interested in the design of a smooth feedback control $u(x)$ that will stabilize the system to any point $x_0 \in \mathcal{E}_d$. We call these points the *design targets* or *design equilibria*:

$$\mathcal{E}_d = \{x_0 \in M \text{ s.t. } x_0 \text{ is a design equilibrium}\}$$

Let

$$\mathcal{E} = \{x_0 \in M \text{ s.t. } f(x_0, u(x_0)) = 0\},$$

the set of equilibria of (4). We assume that $\mathcal{E}$ is finite.

As explained above, when the system evolves on a non-trivial manifold, the Morse inequalities make it unreasonable to expect that there exists a control $u(x)$ that makes the design equilibria the *only* equilibria of the system, i.e. such that $\mathcal{E}_d = \mathcal{E}$. We call these additional equilibria, that are introduced by the non-trivial topology of the space, *ancillary equilibria*:

$$\mathcal{E}_a = \mathcal{E} - \mathcal{E}_d.$$

**Example 2.** *In the case of the 2-cycles problem, there are five distance constraints to satisfy. As was shown in Part I, four formations in the plane, modulo rotation and translation, satisfy the distance constraints; hence, $|\mathcal{E}_d| = 4$. Given $d \in \mathcal{L}$, we call* frameworks *attached to d the frameworks whose edge lengths are given by d. Equivalently, frameworks attached to d are given by letting the reflections $R_1$ and $R_2$ defined in Section 4 or part I act on any framework with edge lengths d. We illustrate these frameworks in Figure 3. Notice that the mirror-symmetric of a formation cannot be obtained via rotations and translations.*

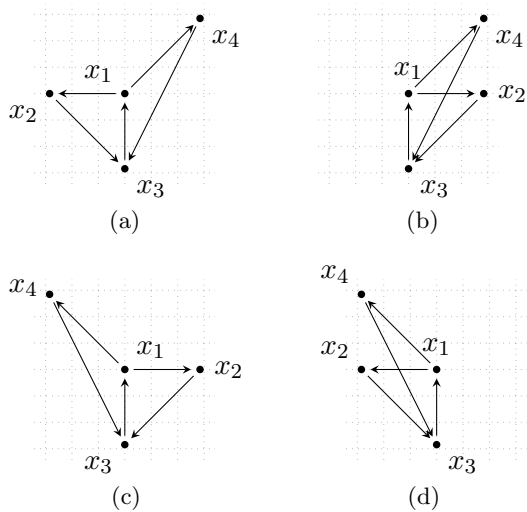

Figure 3: Four formations in the plane that are not equivalent under rotations and translation and that have the same corresponding edge lengths. $(a)$ is the mirror-symmetric of $(c)$ and $(b)$ is the mirror-symmetric of $(d)$.

Let us assume for the time being that the linearization of the system at an equilibrium has no eigenvalues with zero real part. We decompose the set $\mathcal{E}$ into *stable* equilibria, by which we mean equilibria such that *all the eigenvalues* of the linearized system have a negative real part, and *unstable equilibria*, where *at least one eigenvalue* of the linearization has a positive real part. Observe that under this definition, saddle points are considered unstable; this is motivated by practical considerations since the set of initial conditions that result in a system settling at a saddle point is of measure zero[1]. Equivalently, perturbations will move the system away from the $\alpha-$set [GH83] of a saddle point with probability one. In summary:
$$\mathcal{E} = \mathcal{E}_s \cup \mathcal{E}_u$$

---

[1] for e.g. the Lebesgue measure on $M$

where
$$\mathcal{E}_s = \{x_0 \in \mathcal{E} \mid x_0 \text{ is stable}\}$$
and
$$\mathcal{E}_u = \{x_0 \in \mathcal{E} \mid x_0 \text{ is unstable or a saddle}\}.$$

With these notions in mind, we introduce the following definition:

**Definition 1.** *Consider the smooth control system $\dot{x} = f(x, u(x))$ where $x \in M$ and the set $\mathcal{E}$ of equilibria of the system is finite. Let $\mathcal{E}_d \subset M$ be a finite set. We say that $\mathcal{E}_d$ is*

1. *feasible if we can choose a smooth $u(x)$ such that $\mathcal{E}_d \cap \mathcal{E} \neq \emptyset$.*

2. *type-A stable if we can choose a smooth $u(x)$ such that $\mathcal{E}_s \subset \mathcal{E}_d$.*

3. *strongly type-A stable if we can choose a smooth $u(x)$ such that $\mathcal{E}_s = \mathcal{E}_d$.*

*When the set $\mathcal{E}_d$ is clear from the context, we say that the system is type-A stable.*

The set $\mathcal{E}_d$ is feasible if we can choose $u(x)$ such that *at least one* equilibrium of the system is a design target. It is said to be *type-A* stable if the system stabilizes to $\mathcal{E}_d$ with probability one for any randomly chosen initial conditions on $M$. It is *strongly type-A* stable if it is type-A stable and moreover all elements of $\mathcal{E}_d$ are stable equilibria.

The usual notion of global stability is a particular instance of type-A stability; indeed, it corresponds to having $u(x)$ such that $\mathcal{E}_d = \mathcal{E} = \mathcal{E}_s$. Looking at the contrapositive of this definition, a system is *not type-A stable* if there exists a set of initial conditions, which is of strictly positive measure, that lead to an ancillary equilibrium. We observe that type-A stability is a global stability notion; in particular, if one can choose $u$ such that all design equilibria are locally stable, but if this choice forces the appearance of other, undesired equilibria which are also locally stable, the system is not type-A stable. The examples below illustrate these notions.

**Example 3.** *Consider a system*
$$\dot{x} = x(1 - kx^2)$$
*where $k \in \mathbb{R}$ is a feedback parameter to be chosen by the user. We show that any $\mathcal{E}_d \subset (0, \infty)$ is not type-A stable. We first observe that the system has an equilibrium at $0$ and two equilibria at $x = \pm\sqrt{1/k}$ if $k > 0$. The system is thus feasible for any $\mathcal{E}_d \subset \mathbb{R}$. The Jacobian of the system is $1$ at $x = 0$ and $-2$ at $x = \pm\sqrt{1/k}$. For $k > 0$, the above says that*
$$\mathcal{E} = \{0, \pm\sqrt{1/k}\} = \underbrace{\{\sqrt{1/k}\}}_{\mathcal{E}_d} \cup \underbrace{\{0, -\sqrt{1/k}\}}_{\mathcal{E}_a}.$$

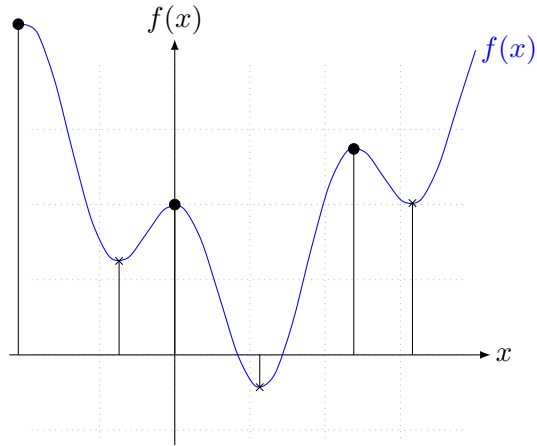

Figure 4: The set $\mathcal{E}_s$ is the set of local minima (marked by ×); the set $\mathcal{E}_u$ is the set of saddles and local maxima (marked by ·). A gradient flow is type-A stable if and only if $\mathcal{E}_d$ contains the set of local minima.

*From the linearization of the system, we have that*

$$\mathcal{E}_s = \{\pm\sqrt{1/k}\} \text{ and } \mathcal{E}_u = \{0\}.$$

*We conclude that $\mathcal{E}_s \not\subseteq \mathcal{E}_d$ and the system is not type-A stable.*

**Example 4** (Gradient flows). *Let $f(x) : \mathbb{R}^n \to \mathbb{R}$ be a smooth function. Assume that a control system is given by the gradient flow of $f$:*

$$\dot{x} = -\frac{\partial f}{\partial x}.$$

*We have that $\mathcal{E}_s$ is the set of local minima of $f$ and $\mathcal{E}_u$ contains all the other points $x_0$ such that $\frac{\partial f}{\partial x}|_{x_0} = 0$. Hence, $\mathcal{E}_d$ is type-A stable if and only if it contains the set of local minima of the gradient flow. We illustrate this in Figure 4.*

## 3 Singularities, genericity and jet spaces

Informally speaking, a property of elements of an arbitrary set is said to be *generic* if it is shared by *almost all* elements of the set. The set in question could be a manifold if one is considering all possible initial conditions of a Cauchy problem, or $\mathcal{C}^2(M)$ if one is considering the class of twice-differentiable functions on $M$. Precisely, we have the definition:

**Definition 2.** *A property $\mathcal{P}$ is* generic *for a topological space $S$ if it is true on an everywhere dense intersection of open sets of $S$.*

Everywhere dense intersections of open sets are sometimes called *residual* sets [AAIS94]. In general, asking for a given property to be generic is a rather strong requirement, and oftentimes it is enough to show that a given property is true on an open set of parameters, initial conditions, etc. We define

**Definition 3.** *An element $u$ of a topological space $S$ satisfies the property $\mathcal{P}$* robustly *if $\mathcal{P}$ is true in a neighborhood of $u$ in $S$.*

If the property is satisfied for a *non-robust $u$*, then it fails to be satisfied under the slightest error in modelling or measurement. We also say that a *property is robust* to refer to the existence of a robust $u$ which satisfies the property.

**Remark 1.** *We emphasize that when we seek a robust control law $u(x)$ for stabilization, we seek a control law such that the equilibrium remains stable under small perturbations in $u(x)$. The equilibrium, however, may move in the state space. For example, assume that the system*
$$\dot{x} = u(x)$$
*has the origin as a stable equilibrium. If for all $g(x)$ in an appropriate set of perturbations, the system*
$$\dot{x} = u(x) + \varepsilon g(x)$$
*has a stable equilibrium at a point $k(\varepsilon)$ near the origin, then the control law $u(x)$ is robust. If, on the contrary, the equilibrium disappears or becomes unstable, then $u(x)$ is not robust.*

Genericity will appear in two guises in this work. The first appearance was in the previous section: we can paraphrase the definition of type-A stability by saying that a system is type-A stable if it will evolve to the set $\mathcal{E}_d$ *generically* with respect to the choice of initial conditions. The other is the subject of this section: we will restrict our search to *robust control laws $u(x)$* for a property $\mathcal{P}$, where $\mathcal{P}$ is either type-A stability or strong type-A stability.

If $\rceil \mathcal{P}$, the opposite of $\mathcal{P}$, is generic, then there is no robust $u$ that satisfies $\mathcal{P}$. Indeed, if $\rceil \mathcal{P}$ is generic, then $\mathcal{P}$ is verified on at most a nowhere dense closed set. In particular, $\mathcal{P}$ is not verified on an open set. We present the main tools to handle genericity and robustness in Appendix B.

# 4 The Logistic equation

We now recall a few definitions from dynamical systems theory. Consider a dynamical system of the form

$$\dot{x} = f_\mu(x) \tag{5}$$

where $x \in M$, an $n-$dimensional manifold, and $\mu \in \mathbb{R}^k$ is a vector of parameters on which the system smoothly depends.

**Definition 4** (Hyperbolic equilibrium and bifurcation value)**.**

1. *An equilibrium $x_0$ is called* hyperbolic *if the eigenvalues of the linearization at $x_0$ have non-zero real-parts.*

2. *A value $\mu_0$ in the parameter space $\mathbb{R}^k$ for which the flow of* (5) *has an eigenvalue with zero real-part is called a* bifurcation value.

Let us put the definition above in context. Assume that $x_0$ is an equilibrium. If it is hyperbolic, the dynamics in a neighborhood of $x_0$ is determined by the number of positive and negative eigenvalues of $\frac{\partial f}{\partial x}|_{x_0}$. Precisely, the Hartman-Grobman theorem [Gro59, Har60] asserts the existence of a homeomorphism that takes the flow in a neighborhood of $x_0$ to the flow of

$$\dot{z} = \begin{bmatrix} I_l & 0 \\ 0 & -I_m \end{bmatrix} z$$

where $I_l$ is the $l \times l$ identity matrix and $l, m$ are the number of positive and negative eigenvalues of $\frac{\partial f}{\partial x}|_{x_0}$, respectively. The image of $[z_1, \ldots, z_l]^T$ under the homeomorphism is called the unstable manifold, and the image of $[z_{l+1}, \ldots, z_{l+k}]^T$ the stable manifold. Hence, under the assumption that an equilibrium is hyperbolic, the behavior of the system in a neighborhood of that equilibrium is entirely determined by its linearization.

At least since the work of Andronov and Pontryagin [AP37], it has been recognized that the behavior of systems around bifurcation points can be far more complex than the usual stable/unstable manifolds decomposition of hyperbolic dynamics [GH83]. As examples of the wide variety of situations can occur nearby a bifurcation value, we mention that periodic orbits may appear (e.g. in the Andronov-Hopf bifurcation) or the dimensions of the stable and unstable manifolds may change (e.g. the transcritical bifurcation).

In general, such singularities of vector fields are not generic since a small perturbation of $f$ will make its Jacobian non-singular (see Corollary 2 below). We will see though that the *distributed nature of the system makes the existence of such singularities generic*. We come back to this point in Section 8.

We show below that the 2-cycles behaves similarly to the logistic equation, which is presented here, in the sense that they both exhibit the same type of singularities or bifurcation. The logistic equation, which is often used to describe systems in which two competing effects—such as supply and demand or predator and prey— are at play, is the one-dimensional ODE given by
$$\dot{x} = x(\mu - x). \tag{6}$$
This equation displays what is called a *transcritical* or *transfer of stability* bifurcation at $\mu = 0$, which we explain here. Observe that it has two equilibria, one at $x = 0$ and one at $x = \mu$, which coalesce when $\mu = 0$. The linearization of the system about $x$ is
$$df = (\mu - x) - x = \mu - 2x.$$
From this linearization, we see that for $\mu > 0$, the equilibrium $x = 0$ is unstable whereas the equilibrium $x = \mu$ is stable. The situation is reversed for $\mu < 0$. We conclude that at the bifurcation value $\mu = 0$, the two equilibria coalesce and *exchange their stability properties*. We depict the above in Figure 5. This figure is to be compared to Figure 14.

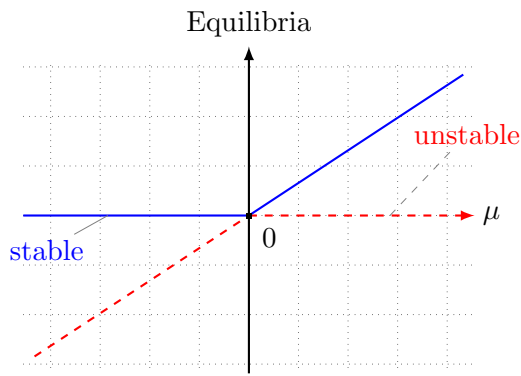

Figure 5: The logistic equation undergoes a transcritical bifurcation when $\mu = 0$. The equilibrium $x = 0$ is stable for $\mu < 0$ and unstable for $\mu > 0$.

The most common approach used to gain some understanding about the behavior of a dynamical system near a non-hyperbolic equilibrium relies on the use of the *center manifold theorem* [GH83]. This theorem asserts the existence of a nonlinear change of coordinates, valid near the equilibrium, where the dynamics can be put in a so-called normal form which is more amenable to analysis. While very useful in general, such an approach is without much hope for success for our objective. Indeed, the change of variables involved in the analysis will depend on the control $u$, and tracking the effect of this dependence through the whole procedure is not feasible for broad classes of control laws. Specifically, in order to show that the 2-cycles behaves similarly to the logistic equation generically for $u$ using a center manifold approach, one would have to exhibit a $u$-dependent nonlinear change of coordinates near singular formations that reduces the 2-cycles to the normal form (6).

In order to sidestep this difficulty, we have recourse to the following result of Sotomayor [Sot73], which characterizes the generic behavior of dynamical systems near non-hyperbolic fixed-points *without recourse* to the center manifold.

**Theorem 1** (Sotomayor). *Let $\dot{x} = f_\mu(x)$ be a system of ODE in $\mathbb{R}^n$ depending on a scalar parameter $\mu$. For $\mu = \mu_0$, assume that the system has an equilibrium $x_0$ satisfying the following conditions:*

1. *$\frac{\partial f_{\mu_0}}{\partial x}|_{x_0}$ has a unique zero eigenvalue with left and right eigenvectors $w$ and $v$ respectively. The other eigenvalues are negative.*

2. *$w^T \frac{\partial f_\mu}{\partial \mu}|_{x_0,\mu_0} v = 0$*

3. *$w^T \frac{\partial^2 f_{\mu_0}}{\partial x^2}|_{x_0}(v,v) \neq 0$ and $w^T \frac{\partial^2 f_\mu}{\partial x \partial \mu}|_{x_0,\mu_0}(v,v) \neq 0$*

*Then the phase portrait is topologically equivalent to the phase portrait of the logistic equation, i.e. we have a transcritical bifurcation about $x_0$ for $\mu = \mu_0$. Thus around $\mu = \mu_0$, there are two arcs of equilibria whose stability properties are exchanged when passing through $\mu_0$. Moreover, the set of equations $\dot{x} = f_\mu(x)$ which satisfy conditions $(1), (2)$ and $(3)$ above is generic in the space of smooth one-parameter families of vector fields with an equilibrium at $x_0$, $\mu_0$ with a zero eigenvalue.*

We will use this theorem to show that the 2-cycles admits a transcritical bifurcation generically and that, as a consequence, $\mathcal{E}_s \nsubseteq \mathcal{E}_d$.

## 5 Algebraic aspects of distributed information

We establish in this section an algebraic framework that reduces questions in decentralized control problems to algebraic ones. Informally speaking, we capture in an algebraic object the range of behaviors that can be achieved given the information flow of a decentralized system. We illustrate this approach in the next section to show that the 2-cycles is generically not strongly type-A stable.

Consider a general decentralized feedback system of the type

$$\begin{cases} \dot{x} &= \sum_i u_i(y_i) g_i(x) = G_u(x) \\ y_i &= h_i(x) \end{cases} \quad (7)$$

where the $g_i(x)$ are smooth vector fields, and the $u_i$ are smooth real-valued functions that belong to a set $\mathcal{U}$ of admissible control laws. The $h_i(x)$ are the observations on the system. What makes the system *decentralized* is the fact that through $h_i$, each control term $u_i$ only has access to partial information about the system. The $h_i$'s define an information flow together with the $u_i$'s. In the case of formation control problems, we have $h_i(x) = \|x\|^2 - d_i$ or $h_i(x) = (e_j, e_k, (x_k - x_i) \cdot (x_j - x_i))$ depending on the number of co-leaders of the agent.

Given a decentralized control system and an objective to attain, there are two notions that need to be formalized in order to know if the objective is feasible:

- What is the range of behaviors that can be achieved by letting $u \in \mathcal{U}$?

- What are the different ways to achieve the objective?

The first notion calls for the description of a reachable set in the space of vector fields, or said differently, a characterization of the possible right-hand sides of Equation (7). The second notion, on which we elaborate below, acknowledges the fact that there often is a *class of system* which satisfy an objective. The class is sometimes explicitly given, e.g. when asking for control laws that will make some eigenvalues negative, or it can be implicitly given, e.g. when requiring that a function vanishes at a point without specifying the order with which it should vanish.

We propose a framework that addresses these two aspects in an algebraic manner. Since the work of Riemann, it has been recognized that there is a duality to be exploited between topological spaces and functions on these spaces. This change of point of view—replacing the study of a space $M$ by the study of functions on that space—has proven remarkably fruitful. Indeed, it has become standard in algebraic geometry, where a variety is studied through sheaves of regular functions, and it is at the core of noncommutative geometry, where the topological space is all but replaced by a noncommutative algebra of functions. This point of view is best summarized by a consequence of the celebrated Gelfand-Neimark theorem asserting that if $M$ is compact and Hausdorff, we can recover $M$ from the knowledge of the ring of smooth complex-valued functions on $M$.

As a consequence, the study of all control problems on $M$ could be reduced, at least in principle, to the study of algebraic objects in spaces of functions on $M$. This approach would rely on translating the vector field/differential description of the system and the objective to a functional description. While this may not bring much to general control systems, we believe it can be very fruitful to answer questions about decentralized control systems. We develop such an approach here. Precisely, we propose a framework that yields algebraic equivalents to the two notions above and reduce the feasibility of an objective to an intersection condition.

We start by reviewing some notions from algebra. A *ring* $R$ is a set with two operations, addition and multiplication, which have the following properties: $R$ is an abelian group for the addition (i.e. the addition is commutative, associative, and there is an additive inverse for every element); the multiplication operation is associative and it has an identity element, but no inverses in general. Moreover, the multiplication and addition operations are distributive. We will work here mostly with the commutative ring $\mathcal{C}^\infty(M)$ of smooth real-valued functions on a manifold $M$ where addition and multiplication are taken pointwise.

Let us illustrate the above point of view on a simple example. Consider the scalar system
$$\begin{cases} \dot{x} &= f(x) + k(x)u(y) = g(x) \\ y &= x^2 \\ u &\in \mathcal{U} \end{cases}.$$

Consider the space of functions

$$\langle g(x)\rangle_{\mathcal{U}} = \left\{g(x) \in \mathcal{C}^\infty(\mathbb{R}) \text{ s.t. } g(x) = f(x) + k(x)u(x^2) \text{ for some } u \in \mathcal{U}\right\}.$$

This space of functions encodes all the possible behaviors of the system. For example, let $p(x)$ be a smooth function such that $p(x_0) = 0$. If $p(x) \in \langle g(x)\rangle_{\mathcal{U}}$, then there exists a feedback such that $x_0$ is an equilibrium for the system. The different ways in which feedback can create an equilibrium at $x_0$ is reflected in the different functions in $\langle g(x)\rangle_{\mathcal{U}}$ that vanish at $x_0$.

We formalize this as follows: recall that an *ideal* $\mathcal{I}$ in a ring $R$ is a subset of $\mathcal{R}$ that is closed under addition and has the property that

$$f \in R \text{ and } h \in \mathcal{I} \Rightarrow fh \in \mathcal{I}.$$

Consider the ideal $I(x_0)$ of functions that vanish at $x_0$. In general, given a subspace $C \subset M$, we define the ideal of functions that vanish on $C$:

$$I(C) = \{g \in \mathcal{C}^\infty(M) \text{ such that } g(x) = 0 \text{ for all } x \in C.\}$$

This is indeed an ideal since $g \in I(C)$ implies $f(x)g(x) = 0$ for all $x \in C, f \in \mathcal{C}^\infty(M)$ and $I(C)$ is closed under addition. A necessary and sufficient condition for the existence of a feedback such that $x_0$ is an equilibrium is thus given by

$$\langle g\rangle_{\mathcal{U}} \cap I(x_0) \neq \emptyset.$$

The above condition effectively replaces the scalar equation $f(x) + k(x)u(x^2) = 0$: *we have thus substituted to a scalar equation in unknown $u$ a problem involving the intersection of function spaces with ideals.* Beyond the fact that an algebraic approach lends itself to concise statements for feasibility, another advantage of translating the problem into algebraic terms is that *there exist well-developed methods to check whether an element belongs to an ideal* for many classes of functions. In case the problem is restricted to *analytic* control systems, there exists a large body of literature on computational methods to handle the objects defined above (e.g. Gröbner Bases) [Sch03].

Our approach is, in broad strokes, the following:

1. Encode the task of the decentralized control problem as a proper ideal in the ring of smooth functions $\mathcal{C}^\infty(M)$. Call this ideal $\mathcal{I}$.

2. Encode the range of behaviors, *with respect to the task at hand*, that the system can achieve via feedback controls $u$ depending on the information flow as a set of functions $\langle F\rangle_{\mathcal{U}}$.

3. Conclude that a necessary condition for the existence of $u$'s that achieve the task is that $\langle F\rangle_{\mathcal{U}} \cap \mathcal{I} \neq \emptyset$.

We derive now a simple necessary criterion for feasibility that will be used in the next section. Let $F \in \mathcal{C}^\infty$ be with the convention that $F(x) \geq 0$ is a necessary condition for the objective to be met. We say that $F$ is *feasible* if we can choose $u$'s such that $F(x) \geq 0$. The situation we have in mind is the one of stability, where we take $F_i(x)$ to be minus the $i^{\text{th}}$ eigenvalue of the linearization at an equilibrium $x$. Hence, $F_i(x) \geq 0$ is a necessary condition for stability.

Assume that at an equilibrium $x$, we have

$$F(x) = p(x)q(x) \tag{8}$$

where $q(x)$ depends on the $u$'s and $p(x)$ is independent of $u$. We define by $V(p)$ the zero set of $p(x)$, i.e.

$$V(p) = \{x \in M \text{ s.t. } p(x) = 0\}$$

Thus $I(V(p))$ is the ideal in $\mathcal{C}^\infty(M)$ of functions that vanish where $p$ vanishes.

Informally speaking, a function $p(x)$ is transversal to a submanifold $C$ of $M$ if the columns of $\frac{\partial}{\partial x} p(x)$ and a basis for the tangent space of $C$ spans the tangent space of $M$. For example, the function $y = x^2$ is not transversal to the submanifold $y = 0$ of $\mathbb{R}^2$ at 0, but the function $y = x$ is. For a formal definition, see Definition 6 in Appendix B. We have the following result:

**Lemma 1.** *Assume that we can factor $F(x)$ as $F(x) = p(x)q(x)$ with $p$ independent of $u$. Let $p(x)$ by such that $V(p)$ is of codimension 1. If $p(x)$ is transversal to $V(p)$, then $F \geq 0$ is feasible only if*

$$\langle q \rangle_\mathcal{U} \cap I(V(p)) \neq \emptyset.$$

*Proof.* Because $p(x)$ is transversal to $V(p)$ at $x$ and $V(p)$ is of codimension 1, $p(x)$ changes sign around $x$. In order for $F(x)$ to remain positive, $q(x)$ needs to vanish at $x$ and thus belong to $I(V(p))$. ∎

Recall that a function is regular at $x_0$ if its Jacobian is of full rank at $x_0$; $x_0$ is also called a regular point for $h$. Let $h : M \to \mathbb{R}^k$ where $M$ is of dimension $n > k$. By the implicit function theorem, the set defined by $h(x) = h(x_0)$ has the structure of a manifold in a neigborhood $V_{x_0}$ of $x_0$. We call this manifold $N_h(x_0)$:

$$N_h(x_0) = \{x \in V_{x_0} \subset M \text{ s.t } h(x) = h(x_0)\}. \tag{9}$$

**Lemma 2.** *Let $C$ be a smooth submanifold of dimension $m$ in an $n$-dimensional manifold $M$. Let $h(x) : M \to \mathbb{R}^k$ and $u(y) : \mathbb{R}^k \to \mathbb{R}$ be smooth functions such that*

$$u(h(x)) \in I(C).$$

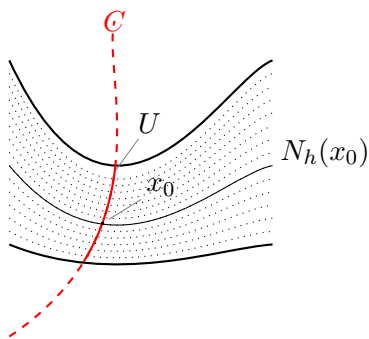

Figure 6: The submanifold $N_h(x_0)$ (the isolevel set of $h$ at the value $h(x_0)$) and $C$ intersect transversally at $x_0$. The dotted curves represent $N_h(c)$ for $c$ in a neighborhood $U \subset C$ of $x_0$. Because $u(h(c))$ vanishes on $U$ and because the value of $u(h(x))$ is constant on all $N_h(c)$ by construction, $u(h(c))$ vanishes on an open neighborhood of $x_0$ in $M$.

*Let $x_0 \in C$ be a regular value for $h$. If $N_h(x_0)$ intersects $C$ transversally at $x_0$, then $u(h(x))$ is zero on an open set in $M$ and $u(y)$ is zero on an open set in $\mathbb{R}^k$. In particular, if we restrict $u$ to be analytic, then $u \equiv 0$..*

*Proof.* The proof is straightforward; in the interest of space, we only sketch it. Without loss of generality, we can take $M = \mathbb{R}^n$, $x_0 = 0$, $h(x_0) = 0$, $C = \mathbb{R}^m$ and $N_h(0) = \mathbb{R}^{n-k}$. Assume first that $C$ and $N_h(0)$ intersect only at 0:

$$C \cap N_h(0) = \{0\};$$

because they intersect transversally, we have that $m = n - k$ . Let $U$ be a neighborhood of 0 in $C$ and let $c \in U$. Because 0 is a regular value of $h$, a neighborhood of 0 in $M$ is contained in $\cup_{c \in U} N_h(c)$. Since $u(h(x)) \in I(C)$ and $c \in U$, then $u(h(c)) = 0$ for $c \in U$. By definition of $N_h$, for all $x \in N_h(c)$, we have $h(x) = h(c)$ and thus $u(h(x))$ vanishes on $\cup_{c \in U} N_h(c)$ which proves the first part.

For the second part, since 0 is a regular value of $h$, $h(x)$ maps a neighborhood of zero in $\mathbb{R}^m$ onto a neighborhood $V$ of 0 in $\mathbb{R}^k$. The previous point thus implies that $u$ vanishes on $V$.

In case $C$ and $N_h$ intersect on a subspace of strictly positive dimension, the argument is the same as the one above but applied to a space $C'$ that is contained in $C$ but whose intersection with $N_h$ is $\{0\}$. We omit the details. The last part uses the fact that the only analytic function that vanishes on an open set is the constant 0. ∎

This lemma is useful to prove that there are no robust elements in an ideal of functions. Indeed, if its conditions are met, by Corollary 2, $u(x)$ is not robust since it vanishes with a zero derivative. The proof that the 2-cycles formation is not robustly strongly type-A stable relies on the same idea as the one developed in the example below, albeit applied to a more complex ideal.

**Example 5.** *Let $f(x,y) = u_1(x)u_2(y)$ and $C_1$ be the subspace of $\mathbb{R}^2$ defined by the equation $x - y = 0$. There are no robust choices of $u_1$ and $u_2$ such that $f(x,y) \in I(C)$.*

*Proof.* Let $h_1(x,y) = x$ and $h_2(x,y) = y$. Since the Jacobians of $h_1$ and $h_2$ are always of full rank, every point in $\mathbb{R}^2$ is a regular value for $h_1$ and $h_2$. Because all the functions involved are continuous, if $f$ is to vanish everywhere on $C$, there exists a connected open set $V$ in $C$ such that either $u_1(h_1(x,y))$ or $u_2(h_2(x,y))$ vanishes. Assume without loss of generality that $u_1(h_1(x)) \in I(V)$ and let $x_0 \in V$. The isolevel sets of $h_1$ are translations of the y-axis in $\mathbb{R}^2$, hence $N_{h_1}$ intersects $C$ transversally everywhere. From Lemma 2, we deduce that $u_1$ is zero in an open set in $\mathbb{R}$ and hence by Corollary 2, it is not robust. ∎

We conclude this section by recording a few simple properties of the operations $V$ and $I$ that we will use in the sequel. We refer the reader to [Sch03] for more details. Assume that $p_1$ and $p_2$, smooth functions on $M$, are such that

$$V(p_1) \subset V(p_2) \Rightarrow (V(p_2)) \subset I(V(p_1)). \tag{10}$$

Moreover, observe that

$$V(p_1 p_2) = V(p_1) \cup V(p_2) \tag{11}$$

since the product $p_1(x)p_2(x)$ of course vanishes if and only if $p_1(x)$ or $p_2(x)$ vanish. Similarly,

$$I(V(p_1) \cup V(p_2)) = I(V(p_1)) \cap I(V(p_2)). \tag{12}$$

# 6 Local stability properties

We show that the 2-cycles formation is not robustly strongly type-A stabilizable or, equivalently, that there does not exist robust feedback laws that will stabilize its four design equilibria. We established in part I a few conditions a feedback control law had to satisfy in order to yield a well-defined formation control system (Definition 7 in part I). We recall them here:

**Definition 7, part I:** *A feedback control law $u_i$ is* compatible *with a formation control problem if*

1. *$u_i(d_j; e_j)$ is such that $u_i(d_j; 0) = 0$ if agent $i$ has one co-leader.*

2. *$u_i(d_j, d_k; e_j, e_k, z_j \cdot z_k)$ is such that $u_i(d_j, d_k; 0, 0, z) = 0$ for all $z$ if agent $i$ has two co-leaders.*

We accordingly define the class of controls $\mathcal{U}$ to be all smooth control laws such that $u_i(d_i; e_i) = 0$ and $u_j(d_i, d_j; e_i, e_j, \cdot) = 0$ for $e_i = e_j = 0$, depending on whether the agent has one or two co-leaders.

Let $x_i \in \mathbb{R}^2$, $i = 1\ldots 4$ represent the positions of the agents. Recall the definition of the vectors

$$\begin{cases} z_1 &= x_2 - x_1 \\ z_2 &= x_3 - x_2 \\ z_3 &= x_1 - x_3 \\ z_4 &= x_3 - x_4 \\ z_5 &= x_4 - x_1 \end{cases} \tag{13}$$

Let $d = (\sqrt{d_1}, \ldots, \sqrt{d_5}) \in \mathcal{L}$ be a target edge length. We use the square roots in order to have $d_i$ enter the expression of the error in edge lengths linearly:

$$e_i = z_i \cdot z_i - d_i.$$

The $e_i$ are lower-bounded by the $d_i$ but not upper-bounded. We define the set of vector fields $\mathcal{F}$ to be the set of admissible vector fields for the 2-cycles formation:

$$\mathcal{F} = \begin{cases} \dot{z}_1 &= u_2(d_2; e_2)z_2 - u_1(d_1, d_5; e_1, e_5, z_1 \cdot z_5)z_1 - u_5(d_1, d_5; e_5, e_1, z_1 \cdot z_5)z_5 \\ \dot{z}_2 &= u_3(d_3; e_3)z_3 - u_2(d_2; e_2)z_2 \\ \dot{z}_3 &= u_1(d_1, d_5; e_1, e_5, z_1 \cdot z_5)z_1 + u_5(d_1, d_5; e_5, e_1, z_1 \cdot z_5)z_5 - u_3(d_3; e_3)z_3 \quad , u_i \in \mathcal{U} \\ \dot{z}_4 &= u_3(d_3; e_3)z_3 - u_4(d_4; e_4)z_4 \\ \dot{z}_5 &= u_4(d_4; e_4)z_4 - u_1(d_1, d_5 e_1, e_5, z_1 \cdot z_5)z_1 - u_5(d_1, d_5; e_5, e_1, z_1 \cdot z_5)z_5 \end{cases} \tag{14}$$

There are four design equilibria for a generic $d \in \mathcal{L}$ (see Figure 3). They are such that $e_i = 0$. We prove the following:

**Theorem 2.** *There is a set of positive measure in $\mathcal{L}$ for which the 2-cycles formation is not robustly strongly type-A stable.*

Precisely, we prove the existence of a set of $d$'s, of strictly positive measure in $\mathcal{L}$, which share the characteristic that at least one of their four design equilibria is unstable for robust $u_i$ in Equation (14).

It has recently been proved [YADF09] that the 2-cycles can be locally stabilized at a given framework in the plane using a relatively simple control law and adjusting some feedback gains. The dynamics used was of the type of Equation 14 with the control law

$$u_i = k_i e_i, i = 1, \ldots, 5$$

where the $k_i$ are real-valued gains used by the agents to stabilize a given framework. It is not guaranteed that such control laws behave well when the system is not near the equilibrium for which the gains have been designed. Theorem 2 states, in contrast, that there are no *globally defined* decentralized feedback law that stabilizes *all* design equilibria locally. A major difference between the two approaches is that the gains in [YADF09] are computed using the *complete knowledge of the framework* whereas we will consider here all possible feedback gains that *respect the decentralized nature of the system*. Our result thus shows that the feedback laws of [YADF09] cannot be computed in a decentralized manner by the agents.

We use for the remainder of this section the following parametrization of $E^4$, the space of planar frameworks with 4 agents: a formation is represented by $\bar{x} = [\bar{x}_2, \bar{x}_3, \bar{x}_4]$ where $\bar{x}_i$ is obtained by translating the formation in order to have $x_1$ at the origin:

$$\bar{x}_2 = x_2 - x_1; \bar{x}_3 = x_3 - x_1; \bar{x}_4 = x_4 - x_1;$$

similarly $\bar{z}_1 = \bar{x}_2 - \bar{x}_1 = \bar{x}_2$, etc. We also use $\bar{x}$ to differentiate these coordinates from the dynamical coordinates $x$: $\bar{x}$ describes an equilibrium or target for the system whereas $x(t)$ represents the position of the system at time $t$.

With the notation of the previous section, the functions

$$h_1(\bar{x}, x) = (d_1, d_5, z_1 \cdot z_5), h_2(\bar{x}) = d_2, h_3(\bar{x}) = d_3, h_4(\bar{x}) = d_4$$

encode the information that is available to agents 1 to 4 respectively *at an equilibrium* to compute gains. Informally speaking, the $e_i$ variables are of no help at design equilibria since they always satisfy $e_i = 0$; they are thus omitted in the $h_i$'s.

We call *local gains* of agents $2, 3$ and $4$ the functions

$$u_{2x}(h_2; 0) = \frac{\partial}{\partial x} u_2(d_2; x)|_{x=0},$$

$u_{3x}(h_3; 0)$ and $u_{4x}(h_4; 0)$. Agent 1 has access to more information and has thus four local gains

$$u_{1x}(h_1(x)) = \frac{\partial}{\partial x} u_1(d_1, d_5; x, y, z_1 \cdot z_5)|_{x,y=0},$$

$u_{1y}(h_1(x))$, $u_{5x}(h_1(x))$, $u_{5y}(h_1(x))$ and $u_{5y}(h_1(x))$.

The solution provided in [YADF09] corresponds to letting each local gain depend on

$$h(\bar{x}) = \bar{x}.$$

We use the notation

$$\begin{aligned}
k_i(d_i) &= u_{ix}(h_i; 0) \text{ for } i = 2, 3, 4 \\
k_{11}(d_1, d_5, z_1 \cdot z_5) &= u_{1x}(h_1(\bar{x}, x)) \\
k_{12}(d_1, d_5, z_1 \cdot z_5) &= u_{1y}(h_1(\bar{x}, x)) \\
k_{51}(d_1, d_5, z_1 \cdot z_5) &= u_{5x}(h_1(\bar{x}, x)) \\
k_{52}(d_1, d_5, z_1 \cdot z_5) &= u_{5y}(h_1(\bar{x}, x))
\end{aligned}$$

for the local gains.

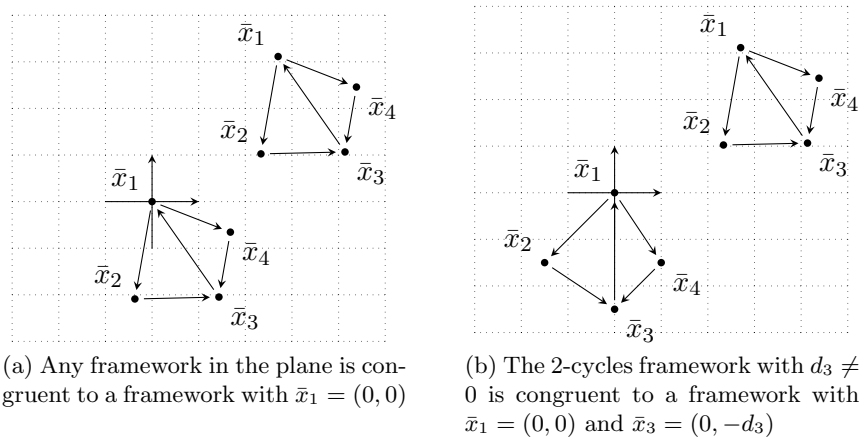

(a) Any framework in the plane is congruent to a framework with $\bar{x}_1 = (0,0)$

(b) The 2-cycles framework with $d_3 \neq 0$ is congruent to a framework with $\bar{x}_1 = (0,0)$ and $\bar{x}_3 = (0,-d_3)$

Figure 7

We prove Theorem 2 in several steps. The state space of minimally rigid formations in the plane is of dimension $2n-3$ (see Laman's theorem, Theorem 1 in Part I). Let $\lambda_1,\ldots,\lambda_{2n-3}$ the eigenvalues of the linearization of the system at an equilibrium $\bar{x}$ and define $F(\bar{x})$ to be minus the determinant of the Jacobian at an equilibrium:

$$F(\bar{x}) = -\prod_{i=1}^{2n-3} \lambda_i.$$

A necessary condition for the equilibrium to be stable is that $F(\bar{x}) \geq 0$. We now prove that $F(\bar{x})$ admits a factorization of the type of Equation (8). Recall the definition of $J$ in Corollary 1 in part I:

$$J = ZA_e Z^T,$$

where $A_e$ is the edge-adjacency matrix of the 2-cycles and

$$Z = \begin{bmatrix} z_1^T & 0 & 0 & \ldots & 0 \\ 0 & z_2^T & 0 & \ldots & 0 \\ 0 & \ldots & & \ddots & \vdots \\ 0 & 0 & \ldots & 0 & z_5^T \end{bmatrix}. \tag{15}$$

Define

$$z_i^\perp = \begin{pmatrix} 0 & 1 \\ -1 & 0 \end{pmatrix} z_i.$$

We have the following formula:

**Lemma 3** (Factorization Lemma). *Let $\bar{x} = [\bar{x}_2, \bar{x}_3, \bar{x}_4]$ represent an equilibrium framework for the 2-cycles. For the dynamics of Equation (14), $F(\bar{x})$ is given by*

$$- F(\bar{x}) = \det(J) = p(\bar{x})q(\bar{x}) \tag{16}$$

*where*

$$q(\bar{x}) = (k_2 k_3 k_4)(k_{11} k_{52} - k_{12} k_{51}) \tag{17}$$

*and*

$$p(\bar{x}) = - \det(A_1)\det(A_2)\det(A_3)\det(A_4)$$

*where the $A_i$ are $2 \times 2$ matrices given by*

$$A_1 = \begin{bmatrix} | & | \\ \bar{z}_1 & \bar{z}_3 \\ | & | \end{bmatrix}, A_2 = \begin{bmatrix} | & | \\ \bar{z}_1 & \bar{z}_5 \\ | & | \end{bmatrix}, A_3 = \begin{bmatrix} | & | \\ \bar{z}_3 & \bar{z}_4 \\ | & | \end{bmatrix}, A_4 = \begin{bmatrix} \bar{z}_1 \cdot \bar{z}_3^\perp & \bar{z}_4 \cdot \bar{z}_3^\perp \\ \bar{z}_2 \cdot \bar{z}_2 & \bar{z}_4 \cdot \bar{z}_4 \end{bmatrix}. \tag{18}$$

*Proof.* We write
$$p_i(\bar{x}) = \det(A_i).$$

The first equality in Equation (16) is a consequence of Corollary 1 in part I. We leave the details to the appendix. ∎

We have thus expressed the determinant of $J$ as the product of a *purely geometric* term, that is $p(\bar{x})$, and a dynamical term, that is $q(\bar{x})$, that depends on the local gains $k_i$. The Factorization lemma is the basis of the proof of Theorem 2 below, whose idea we illustrate in Figure 8.

*Theorem 2.* Recall that $\bar{x} = [\bar{x}_{21}, \bar{x}_{22}, \bar{x}_{31}, \bar{x}_{32}, \bar{x}_{41}, \bar{x}_{42}]$ describes an equilibrium framework as in Figure 7a. Let $F(\bar{x}) = -\prod_{i=1}^{5} \lambda_i$ where the $\lambda_i$'s are the eigenvalues of the linearization of the 2-cycles at an equilibrium framework $\bar{x}$. The system is strongly type-A stable for $d \in \mathcal{L}$ if there exists $k_i$'s, which depend on $d$ via the functions $h_i$, such that $F(\bar{x}) \geq 0$. From Lemma 3, we know that
$$F(\bar{x}) = p(\bar{x})q(\bar{x})$$
where $q$ depends on the choice of feedback gains $k_i$ and $p(\bar{x})$ does not.

We proceed by showing that all the elements of $\langle q(\bar{x}) \rangle_\mathcal{U} \cap I(p)$ are *not robust*. We conclude, using Lemma 1, that there is a set of positive measure of frameworks for which $F(\bar{x}) < 0$. The system is thus not robustly strongly type-A stable for the edge-lengths $d$ of these frameworks.

Let $p_1(\bar{x}) = \det(A_1)$. From Equations (11) and (12), we have that $I(V(p)) \subset I(V(p_1))$. Hence, it is enough to show that elements in $\langle q(x) \rangle_\mathcal{U} \cap I(V(p_1))$ are not robust.

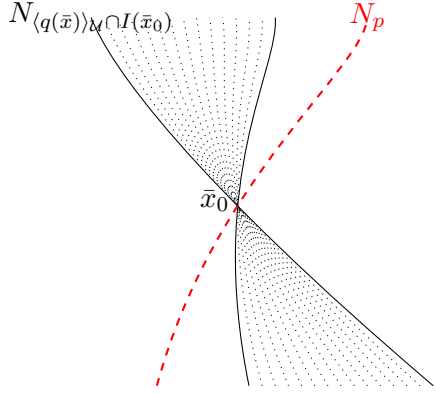

Figure 8: The zero-set of $p(x)$ is depicted by the dashed curve. The zero sets of functions in $\langle q(\bar{x})\rangle_{\mathcal{U}}$ that go through a generic point $\bar{x}_0$ with $p(\bar{x}_0) = 0$ are represented by the dotted curves. A necessary condition for strong type-A stability is that $F(\bar{x}) \geq 0$. Because the no robust function in $\langle q(x)\rangle_{\mathcal{U}}$ has a zero-set that contains the zero-set of $p(\bar{x})$, when $p(\bar{x})$ crosses zero and changes sign, so does $F(\bar{x})$ and the system is thus not strongly type-A stable.

In the coordinate system used, the set $V(p_1)$ is given by

$$V(p_1) = \{\bar{x} = [\bar{x}_2, \bar{x}_3, \bar{x}_4] \text{ s.t. } \det(\bar{x}_2, \bar{x}_3) = 0\}$$

where $\det(x_i, \bar{x}_j)$ is the determinant of the two by two matrix with columns $\bar{x}_i$ and $\bar{x}_j$. This corresponds to frameworks where $x_1, x_2$ and $x_3$ are aligned. Except at the frameworks in $V(p_1)$ where $\det(A_i) = 0$ for $i = 2, 3$ or $4$, $p(\bar{x})$ is transversal to $V(p_1)$. Hence, it is transversal except on a set of codimension at least one in $V(p_1)$, which is necessarily of zero measure by Sard's theorem [War83].

We have that

$$\langle q(x)\rangle_{\mathcal{U}} = \{k_2(h_2)k_3(h_3)k_4(h_4)\left(k_{11}k_{52} - k_{12}k_{51}\right)(h_1)\},$$

where we recall that $h_2(\bar{x}) = d_2(\bar{x})$, $h_3(\bar{x}) = d_3(\bar{x})$, $h_4(\bar{x}) = d_4(\bar{x})$ and $h_1(\bar{x}, x) = (d_1(\bar{x}), d_5(\bar{x}), x_2 \cdot x_4)$.

Let $\bar{x}_0$ be a generic framework in $V(p_1)$. It is possible to perturb $\bar{x}_0$ so that $\bar{x}_0 + \delta\bar{x}$ is not in $V(p_1)$ while $h_2$ is constant—take for example $\delta\bar{x}$ to be, for $\varepsilon$ small,

$$\delta\bar{x}_2 = [0, \varepsilon(\bar{x}_2 - \bar{x}_3)^{\perp}, 0, 0] :$$

this perturbation rotates $\bar{x}_2$ about $\bar{x}_3$. It thus clearly leaves $h_2$ and $h_4$ constant but not $\det(x_2, \bar{x}_3)$. We illustrate it in Figure 9. Because $V(p_1)$ is of codimension 1, we conclude that the tangent space of $N_{h_2}$ (resp. $N_{h_4}$) together with the tangent space of $V(p_1)$ spans the tangent space of $E^4$ at $\bar{x}_0$. Hence, $N_{h_2}$ (resp. $N_{h_4}$) intersect $V(p_1)$ transversally. The same holds for $N_{h_1}$ and $N_{h_3}$, with the perturbation

$$\delta\bar{x}_1 = \left[0, 0, \varepsilon\bar{x}_3^\perp, 0\right],$$

which corresponds to a rotation of $\bar{x}_3$ about $\bar{x}_1$.

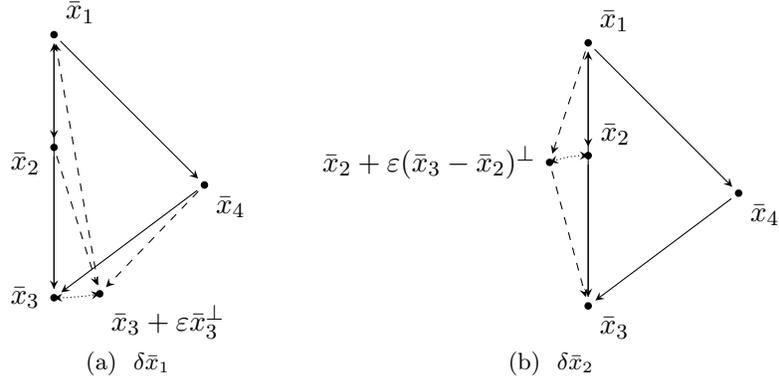

Figure 9: The perturbation $\delta\bar{x}_1$ rotates $\bar{x}_3$ about $\bar{x}_1$ and leaves $d_1, d_3, d_5$ and $z_1 \cdot z_5$ constant but does not leave $x_2, x_3$ and $x_4$ aligned. Similarly, the perturbation $\delta\bar{x}_2$ rotates $\bar{x}_2$ about $\bar{x}_3$ and leaves $d_2$ and $d_4$ constant but does not leave $x_2, x_3$ and $x_4$ aligned.

For $q(\bar{x})$ to vanish on $V(p_1)$, at least one of its factors has to vanish on an open set in $V(p_1)$. Assume that $k_2(d_2)$ vanishes on an open set $U \subset V(p_1)$. By Lemma 2, we know that $k_2(d_2)$ is then zero on a connected open set and hence is not robust by Corollary 2. The same argument applies to the other factors in $q(\bar{x})$. ∎

Informally speaking, the proof above relies on the fact that no agent can decide whether $x_1, x_2$ and $x_3$ are aligned, and hence no agent can take the appropriate action to counter the sign change of $p(x)$. We could have similarly shown that $\langle q(x) \rangle_{\mathcal{U}} \cap I(V(p_3))$ consists of non-robust elements. The same is *not true*, however, for $V(p_2)$. Indeed, $V(p_2)$ corresponds to frameworks such that $x_1, x_2$ and $x_4$ are aligned, and such frameworks can be observed by agent 1. In that case, the argument above fails when trying to prove the transversality of $N_{h_1}$ and $V(p_2)$.

## 6.1 Unstable equilibria

We illustrate in this section Theorem 2 by exhibiting a set of edge lengths in $\mathcal{L}$ that has at least one unstable design equilibrium. We will do so by finding frameworks at which the sign of $p(x)$ changes while the sign of $q(x)$ remains the same.

We have seen in Section 4.5, Part I that there is an action of the group $\mathbb{Z}_2 \times \mathbb{Z}_2$ whose orbit on a given framework leaves the vector $d$ unchanged. Notice that given a vector $d \in \mathcal{L}$, $\det(J)$ is invariant under the mirror symmetry $(1,1)$ (or $R_1 R_2$) but is *not* invariant under the *reflection* symmetries $(1,0)$ and $(0,1)$ (or $R_1$ and $R_2$). The group of transformations given by the reflections $R_1$ and $R_2$ is not rich enough to illustrate Theorem 2. We thus define here another action of $\mathbb{Z}_2 \times \mathbb{Z}_2$ that combines $R_1$ and a reflection defined below.

Frameworks such that $d_3 = 0$ are degenerate, in the sense that $\det(J) = 0$, and form a subset of codimension 1 in $\mathcal{L}$. We thus assume that $d_3 > 0$. Without loss of generality, we can consider the coordinate system depicted in Figure 7b: it consists of aligning $z_3$ with the vertical axis using rotational invariance. In this coordinate system, $\bar{x}_1 = [0,0]$ and $\bar{x}_3 = [0, -d_3]$ and a framework is thus described by $\bar{x}_2, \bar{x}_4$ and $d_3$.

We define $\mathcal{L}_c$ to be the subset of $\mathcal{L}$ corresponding to edge lengths such that at least one corresponding framework has $\bar{x}_1$ in the convex hull of $\bar{x}_2, \bar{x}_3, \bar{x}_4$. With the notation of Section 2 in Part I, we have:

**Definition 5.** *We let $\mathcal{L}_c \subset \mathcal{L}$ be the set of edge lengths $d$ such that for any framework $x_o$ corresponding to $d$, at least one of the frameworks $\bar{x}_o, R_1(\bar{x}_o), R_2(\bar{x}_o), R_1 R_2(\bar{x}_o)$ has $\bar{x}_1$ in the convex hull of $\bar{x}_2, \bar{x}_3, \bar{x}_4$.*

This condition is invariant under the $SE(2)$ action of the framework and $\mathcal{L}_c$ is thus well-defined. See Figure 10 for an example of element of $\mathcal{L}_c$.

**Lemma 4.** *If $d \in \mathcal{L}_c$ then $p(\bar{x}) > 0$ for the framework such that $\bar{x}_1$ is in the convex hull of $\bar{x}_2, \bar{x}_3, \bar{x}_4$.*

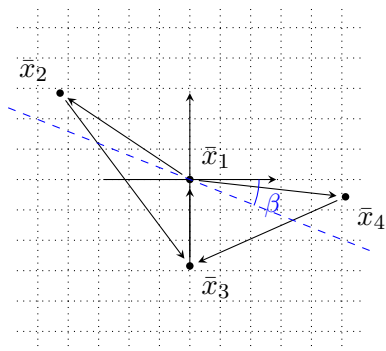

Figure 10

*Proof.* For $\bar{x}_1$ to be in the convex hull of $\bar{x}_2, \bar{x}_3$ and $\bar{x}_4$, $x_{21}$ and $x_{41}$ need to be of opposite signs as illustrated in Figure 10. Without loss of generality, we assume that $\bar{x}_{21} < 0$ and $\bar{x}_{41} > 0$.

The convexity hypothesis implies the existence of $\beta \in \mathbb{R}$, such that

$$\begin{cases} \beta \bar{x}_{21} & < & \bar{x}_{22} \\ \beta \bar{x}_{41} & < & \bar{x}_{42} \end{cases} \tag{19}$$

In these coordinates we have:

$$\bar{z}_1 = \bar{x}_2, \bar{z}_2 = [0, -d_3] - \bar{x}_2, \bar{z}_3 = [0, d_3], \bar{z}_4 = [0, -1] - \bar{x}_4, \bar{z}_5 = \bar{x}_4.$$

We look at the terms $p_1$, $p_2$, $p_3$ and $p_4$ in these coordinates:

1. $p_1 p_3 = d_3^2 \bar{x}_{21} \bar{x}_{41} < 0$:
   This is clear since $\bar{x}_{21} < 0$ and $\bar{x}_{41} > 0$.

2. $p_2 < 0$:
   Recall that $p_2 = \det(A_2)$, where $A_2$ is given in Equation (18). The sign of $p_2$ is thus positive if the vectors $\bar{x}_2$, $\bar{x}_4$ are positively oriented and negative otherwise. Under the above assumption, they are negatively oriented.

3. $p_4 > 0$:
   In the coordinates used, $p_4 = d_3(\bar{x}_{41}d_2^2 - \bar{x}_{21}d_4^2)$. Since $\bar{x}_{21} < 0$ and $\bar{x}_{41} > 0$, we have that $p_4 > 0$.

Under the opposite assumption, namely that $\bar{x}_{21} > 0$ and $\bar{x}_{41} < 0$, the sign of $p_1$ is unchanged and both the signs of $p_2$ and $p_3$ are changed. ∎

We established in part I that $u_z(d_1, d_5; 0, 0, z) = 0$, but we do not have in general that $u_x(d_1, d_5; 0, 0, z) = 0$. We consequently let $h_1$ depend on $d_1, d_5$ and the inner product $x_2 \cdot x_4$, which is equal to $\bar{x}_2 \cdot \bar{x}_4$ at a design equilibrium. To address this dependence, we need to introduce a partial reflection of the framework that keeps this angle constant. We use the same notation as in Section 4 of Part I and define

$$R_3(\bar{x}_1, \bar{x}_2, \bar{x}_3, \bar{x}_4) = (\bar{x}_1, \mathcal{R}^{\bar{z}_3^\perp} \bar{x}_2, \bar{x}_3, \mathcal{R}^{\bar{z}_3^\perp} \bar{x}_4).$$

Observe that $d_1, d_2, d_3$ and $\bar{x}_2 \cdot \bar{x}_4$ are preserved under $R_3$ but $d_2$ and $d_4$, in general, are not. We illustrate $R_3$ in Figure 11. Recall the symmetry $R_1$ which sends $\bar{x}_2$ to its mirror symmetric along $\bar{z}_3$. We have that

$$R_1 R_3 = R_3 R_1, R_1^2 = Id, R_3^2 = Id$$

and thus $R_1$ and $R_3$ generate an action of $\mathbb{Z}_2 \times \mathbb{Z}_2$ on a framework. We show the orbit of this group on a framework in Figure 12. We have the following result:

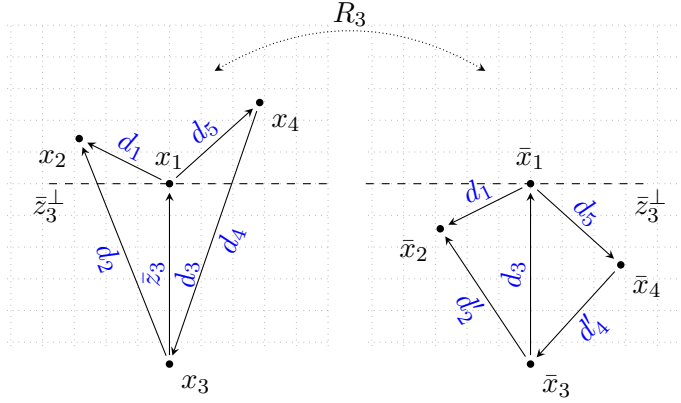

Figure 11: The pseudo-reflection $R_3$ leaves $\bar{x}_1$ and $\bar{x}_3$ fixed and sends $\bar{x}_2$ and $\bar{x}_4$ to their mirror symmetric along $\bar{z}_3^\perp$.

**Lemma 5.** *Let $\bar{x}_o$ be a framework with edge lengths given by $d \in \mathcal{L}$. There is a set of positive measure in $\mathcal{L}$ such that*

$$p(\bar{x}_o) > 0, p(R_1(\bar{x}_0)) < 0, p(R_3(\bar{x}_o)) < 0 \text{ and } p(R_1 R_3(\bar{x}_o)) < 0.$$

*Proof.* Since $p(\bar{x})$ is smooth, it is enough to exhibit a $d$ which satisfies the above relations. Take, e.g., $d = [1, \sqrt{2}, 1, \sqrt{5}, \sqrt{2}]$. A straightforward computation yields the result. ∎

We now prove the following:

**Lemma 6.** *For configurations of the type described in Lemma 5, there are no $u_i \in \mathcal{U}$ such that all design equilibria are exponentially stable.*

*Proof.* Consider a configuration of the type of Lemma 5. Write

$$f(\bar{x}) = k_{11}(h_1(\bar{x}))k_{52}(h_1(\bar{x})) - k_{12}(h_1(\bar{x}))k_{51}(h_1(\bar{x}))$$

and

$$g(\bar{x}) = k_2(h_2(\bar{x}))k_3(h_3(\bar{x}))k_4(h_4(\bar{x})).$$

Recalling that $R_1$ keeps all the $d_i$ fixed but not $\bar{x}_2 \cdot \bar{x}_4$ and that $R_3$ keeps $d_1, d_5, \bar{x}_2 \cdot \bar{x}_4$ fixed, we have

$$\begin{cases} f(R_3(\bar{x})) &= f(\bar{x}) \\ g(R_1(\bar{x})) &= g(\bar{x}) \end{cases}. \tag{20}$$

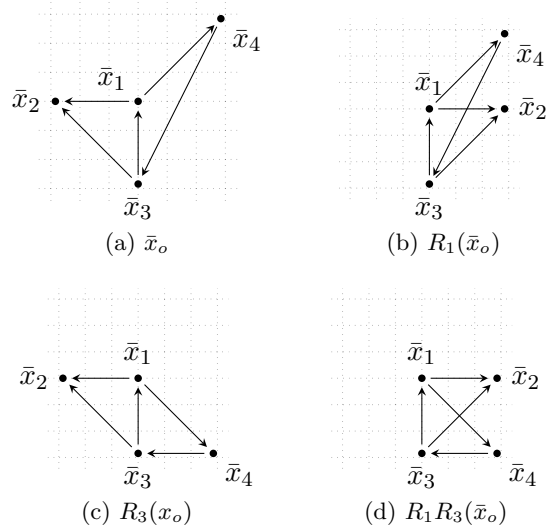

Figure 12: In (a), we represent the framework $\bar{x}_o$ with $d = [1, \sqrt{2}, 1, \sqrt{5}, \sqrt{2}]$. The other figures illustrate the action of $R_1, R_3$ on $\bar{x}_o$.

We gather the conclusions of Lemma 5 and Equation (20) in the following table, where the first two columns and the last two columns correspond to frameworks with the same edge lengths :

|  | $\bar{x}_o$ | $R_1(\bar{x}_o)$ | $R_3(\bar{x}_o)$ | $R_1 R_3(\bar{x}_o)$ |
|---|---|---|---|---|
| $f(\bar{x})$ | $f(\bar{x}_o)$ | $f(R_1(\bar{x}_o))$ | $f(\bar{x}_o)$ | $f(R_1(\bar{x}_o))$ |
| $g(\bar{x})$ | $g(\bar{x}_o)$ | $g(\bar{x}_o)$ | $g(R_3(\bar{x}_o))$ | $g(R_3(\bar{x}_o))$ |
| $p(\bar{x})$ | + | - | - | - |
| product | $f(\bar{x}_o)g(\bar{x}_o)$ | $-f(R_1(\bar{x}_o))g(\bar{x}_o)$ | $-f(\bar{x}_o)g(R_3(\bar{x}_o))$ | $-f(R_1(\bar{x}_o))g(R_3(\bar{x}_o))$ |

The plus and minus signs in the third row represent the sign of $p(\bar{x})$ at the given framework. The product of the elements of each column are gathered in the last row. Hence, the framework $\bar{x}_o$ is exponentially stable only if $f(\bar{x}_0)g(\bar{x}_o)$ is negative and similarly for the other frameworks.

We claim that there does not exist $f$ and $g$, hence $u_i$'s, such that all the products in the last row are negative. Indeed, the first column requires that either $f(\bar{x}_o)$ or $g(\bar{x}_o)$ be negative. Assume that $f(\bar{x}_o) < 0$. Then, starting from the first column, we have $g(\bar{x}_o) > 0$, which implies using the second column that $f(R_1(\bar{x}_o)) > 0$, which in turns yields $g(R_3(\bar{x}_o)) > 0$ using the fourth column. But the third column then requires $f(\bar{x}_o) > 0$, which is a contradiction. Assuming $f(\bar{x}_o) > 0$ yields a similar contradiction. Hence, we conclude that either $\bar{x}_o, R_1(\bar{x}_o)$ or $R_3(\bar{x}_o), R_1R_3(\bar{x}_o)$, or both pairs, are such that $\det(J) > 0$ for at least one framework, which concludes the proof. ∎

## 7 Singularities, transfer of stability and the appearance of stable ancillary equilibria

We prove in this section that the 2-cycles formation is not robustly type-A stable for a set of $d$'s in $\mathcal{L}$ of positive measure. We proceed by showing that the existence of a *stable ancillary equilibria is generic in $\mathcal{F}$*.

### 7.1 Singular formations for $n = 4$ agents

Consider a minimally rigid formation with 4 vertices. Owing to Laman's theorem, it has $2n - 3 = 5$ edges. We single out a particular type of formations which, even though they are infinitesimally rigid, are degenerate in a sense we describe below. Let us denote the edge lengths for the formation by $(\sqrt{d_1}, \ldots, \sqrt{d_5}) \subset \mathcal{L}$ where we recall that $\mathcal{L}$ is the set of feasible edge lengths and we use the square roots so that the $d_i$ enter the $e_i$ linearly. Recall that $R_1$ and $R_2$, defined in Section 4.5 in Part I, generate an action of $\mathbb{Z}_2 \times \mathbb{Z}_2$ on minimally rigid frameworks with 4 agents. The orbit of this group on a framework with edge lengths $d$ gives all frameworks with edge lengths $d$.

We define $\mathcal{S}$ to be set of edge lengths with at least one attached framework such that $z_1$ is parallel to $z_5$, with the notation of Figure 13:

$$\mathcal{S} = \{d \in \mathcal{L} \text{ s.t. } z_1 \text{ parallel } z_5 \text{ for one framework with edge lengths } d.\}.$$

We define $\mathcal{S}_0 = \mathcal{S} \cap \mathcal{L}_0$. We have the following result:

**Lemma 7.** *The following properties of $\mathcal{S}$ hold:*

1. *$\mathcal{S}$ is of codimension one in $\mathcal{L}$*

2. *The formations corresponding to edge lengths in $\mathcal{S}_0$ are infinitesimally rigid.*

*Proof.* For the first part, observe that we can parametrize $\mathcal{S}$ by choosing a feasible $d_1, d_2, d_3$ yielding a triangle $x_1, x_2, x_3$ and one additional parameter giving the signed length of $z_5$, with the sign referring to $z_5$ going in the same direction as $z_1$ or the opposite direction. We see that we need 4 parameters to describe a formation in $\mathcal{S}$ and hence it is of codimension one.

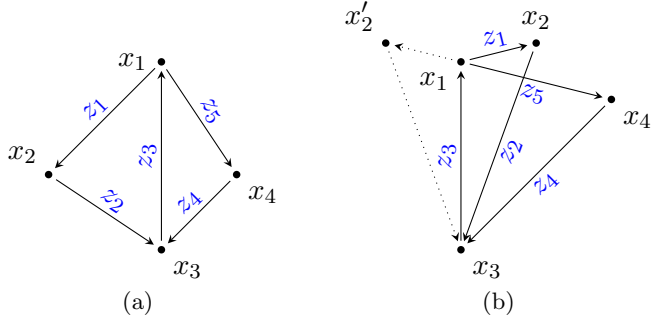

Figure 13: The formation in (a) is such that $(\|z_1\|,\ldots,\|z_5\|) \notin \mathcal{S}$, whereas $(\|z_1\|,\ldots,\|z_5\|) \in \mathcal{S}$ for the framework in (b) since applying $R_1$ to (b) yields a framework with $z_1$ and $z_5$ aligned (dotted lines).

For the second part, we recall that the rigidity matrix of the 2-cycles is given by

$$R = \begin{bmatrix} z_1^T & -z_1^T & 0 & 0 \\ 0 & z_2^T & -z_2^T & 0 \\ -z_3^T & 0 & z_3^T & 0 \\ 0 & 0 & z_4^T & -z_4^T \\ z_5^T & 0 & 0 & -z_5^T \end{bmatrix} \qquad (21)$$

where the $z_i$ are defined in Equation (1). The rigidity matrix $R$ of a framework can be expressed as

$$R = ZA_m^{(2)},$$

where we recall that $A_m$ is the mixed adjacency matrix, $A_m^{(2)} = A_m \otimes I$, $I$ being the $2 \times 2$ identity matrix and $Z$ is as in Equation 15. According to Lemma 4 in part I, the rigidity matrix of the 2-cycles is of full rank when $d \in \mathcal{L}_0$. ∎

## 7.2 Transfer of stability or $\mathcal{E}_s \nsubseteq \mathcal{E}_d$

The set of design equilibria $\mathcal{E}_d$ for the $2-cycles$ is of cardinality 4, up to rigid transformations, since there are four frameworks in the plane for which $e_i = 0$; they are depicted in Figure 3.

Even though the set $\mathcal{E}_a$ of ancillary equilibria depends in general on the choice of feedbacks $u_i \in \mathcal{U}$, some configurations belong to $\mathcal{E}_a$ for all elements of $\mathcal{U}$. We exhibit some of these configurations:

**Proposition 1.** *The set $\mathcal{E}$ contains, in addition to the equilibria in $\mathcal{E}_d$, the frameworks characterized by*

1. $z_i = 0$ for all $i$, which corresponds to having all the agents superposed.

2. all $z_i$ are aligned, which corresponds to having all agents on the same one-dimensional subspace in $\mathbb{R}^2$. These frameworks form a three dimensional invariant subspace of the dynamics.

3. $e_2 = e_3 = e_4 = 0$, $z_1$ and $z_5$ are aligned and so that

$$u_1(e_1, e_5, z_1 \cdot z_5)\|z_1\| = \pm u_5(e_1, e_5, z_1 \cdot z_5)\|z_5\|,$$

where the sign depends on whether $z_1$ and $z_5$ point in the same or opposite directions.

This result is straightforward from an inspection of Equation (14). frameworks of type 2 above are non-infinitesimally rigid and they define an invariant submanifold of the dynamics. An application of the Morse inequalities shows that the 2-cycles dynamics restricted to this invariant submanifold contains at least one equilibrium that is stable. One can nevertheless choose $u$ such that the system will not flow towards this submanifold, and hence this equilibrium can be made unstable for the dynamics in $E^4$. We do not expand on these two ancillary equilibria further here.

The set $\mathcal{S}$ can be alternatively defined as containing the $d$'s for which $\mathcal{E}_d$ contains frameworks of type 3 above. We will show that if one perturbs the value of $d$ around $\mathcal{S}$, then the framework of type 3 will persist in $\mathcal{E}_a$ and moreover keep the stability properties it enjoyed when it was part of $\mathcal{E}_d$.

We now prove that the system of Equation 14 is equivalent to the logistic equation around frameworks in $\mathcal{S}_0$. Denote by $F(z)$ the vector field on the right-hand side of Equation 14, i.e.

$$F(z) = \begin{bmatrix} u_2 z_2 - u_1 z_1 - u_5 z_5 \\ u_3 z_3 - u_2 z_2 \\ u_1 z_1 + u_5 z_5 - u_3 z_3 \\ u_3 z_3 - u_4 z_4 \\ u_4 z_4 - u_1 z_1 - u_5 z_5 \end{bmatrix}, \tag{22}$$

where $u_i$ is $u_i(d_i; e_i)$ for $i = 2, 3, 4$ and $u_i(d_1, d_5; e_1, e_5, z_1 \cdot z_5)$ for $i = 1, 5$. We have the following result

**Theorem 3.** *Consider the set of vector fields*

$$\mathcal{F} = \{F(z) | u_i \text{ are compatible with the 2-cycles}\}$$

*where $F$ is as in Equation (22). A transcritical bifurcation at $\mathcal{S}_0$ is generic for $\mathcal{F}$.*

As a corollary, we will prove

**Theorem 4.** *The 2-cycles formation is not robustly type-A stable for a set of positive measure of parameters $d \in \mathcal{L}$.*

In the interest of clarity, we first prove Theorem 3 for a smaller class of control laws and then show in Appendix D how to extend the proof to the general case. We denote by $F_r(z)$ the vector field

$$F_r(z) = \begin{bmatrix} u(e_2)z_2 - u(e_1)z_1 - u(e_5)z_5 \\ u(e_3)z_3 - u(e_2)z_2 \\ u(e_1)z_1 + u(e_5)z_5 - u(e_3)z_3 \\ u(e_3)z_3 - u(e_4)z_4 \\ u(e_4)z_4 - u(e_1)z_5 - u(e_1)z_1 \end{bmatrix}. \tag{23}$$

This class of vector fields is restricted in the sense that it depends on the $d_i$'s only through the $e_i$'s, and each agents uses the same feedback law. The compatibility condition for $u$ then becomes $u(0) = 0$. We prove

THEOREM 3 BIS. *Consider the set of vector fields*

$$\mathcal{F}_r = \{F_r(z) | u(0) = 0\}$$

*where $F$ is as in Equation* (23). *A transcritical bifurcation at $\mathcal{S}_0$ is generic for $\mathcal{F}_r$.*

Since the relations $z_1 + z_2 + z_3 = z_3 + z_4 + z_5 = 0$ are always satisfied, the $z$ variables used here give a redundant description of the 2-cycles. A set of coordinates that describes the dynamics of the 2-cycles in a non-redundant manner is given by, for example, the error variables $e_i$. These were in fact used in the analysis of the triangular formation in [AYDM07]. In the case of the 2-cycles, since for each value of the $e_i$'s there are generically four possible frameworks with different dynamics, four sets of equations are needed to describe the system in these variables. Explicitly, the $e_i$'s are invariant under $R_1$ and $R_2$, but the dynamics of, say, $x_1$ is not. This need for several set of equations is not particular to the $e_i$, but is a common characteristic of any non-redundant set of variables for the 2-cycles as the LS-category of the state-space of the system is strictly larger than one. See Section 3.3 in Part I. To avoid having to handle several sets of equations, we work with the $z$ variables.

We prove Theorem 3 BIS in several steps.

**Proposition 2.** *Let $d \in \mathcal{S}_0$. There is a non-zero vector $w \in \mathbb{R}^{10}$ such that $w^T \frac{\partial F_r}{\partial z}|_{e_i=0,d} = w^T \frac{\partial F_r}{\partial d}|_d = 0$ for at least one framework attached to $d$.*

We will prove Proposition 2 by relying on some technical lemmas. We use the notation

$$Z_i = z_i z_i^T, Z_i \in \mathbb{R}^{2 \times 2}$$

We evaluated in Proposition 4 (Part I) the linearization of the dynamics of a general formation at a design equilibrium. In the case of the 2-cycles with restricted dynamics, this yields:

$$\frac{\partial F_r}{\partial z}|_{e_i=0} = 2u'(0) \begin{bmatrix} -Z_1 & Z_2 & 0 & 0 & -Z_5 \\ 0 & -Z_2 & Z_3 & 0 & 0 \\ Z_1 & 0 & -Z_3 & 0 & Z_5 \\ 0 & 0 & Z_3 & -Z_4 & 0 \\ -Z_1 & 0 & 0 & Z_4 & -Z_5 \end{bmatrix}. \tag{24}$$

We evaluate the Jacobian of $F_r$ with respect to the vector of edge lengths $d$:

**Lemma 8.** *The Jacobian of $F$ with respect to the parameters $d$ at a design equilibrium is given by*

$$\frac{\partial F_r}{\partial d}|_{e_i=0,d} = -u'(0) \begin{bmatrix} -z_1 & z_2 & 0 & 0 & -z_5 \\ 0 & -z_2 & z_3 & 0 & 0 \\ z_1 & 0 & -z_3 & 0 & z_5 \\ 0 & 0 & z_3 & -z_4 & 0 \\ -z_1 & 0 & 0 & z_4 & -z_5 \end{bmatrix}. \tag{25}$$

*Proof.* We have that $\frac{\partial F_1}{\partial d_1} = -u'_1(e_1)z_1$ and similar expressions for the other entries. ∎

**Lemma 9.** *At a framework attached to $d \in \mathcal{L}_0$, $w$ is a left eigenvector of $\frac{\partial F_r}{\partial z}|_{e_i=0}$ with eigenvalue $0$ if and only if $w^T \frac{\partial F_r}{\partial d} = 0$*

*Proof.* It is easily verified that

$$\frac{\partial F_r}{\partial d} Z = \frac{\partial F_r}{\partial z}, \tag{26}$$

where $Z$ is given in Equation 15.

Since $Z$ is rank deficient only if $z_i = 0$, for any $i$, then $Z$ is of full rank for frameworks attached to elements of $\mathcal{L}_0$. Hence, $v^T Z = 0 \Leftrightarrow v = 0$. We conclude, using Equation 26, that $w^T \frac{\partial F_r}{\partial z} = 0$ if and only if $w^T \frac{\partial F_r}{\partial d} = 0$. ∎

We use Corollary 1 from Part I to obtain the eigenvalues of the linearization of the 2-cycles in the $z$ variables. With the notation of Section 5.3 in part I, we have $Z' = Z$ in the case of the restricted dynamics. We thus conclude that the eigenvalues of (24) are 0 with multiplicity five and the eigenvalues of

$$J = u'(0) Z A_e^{(2)} Z^T \tag{27}$$

where we recall that $A_e$ is the edge-adjacency matrix of the 2-cycles. We deduce from the above equation the following corollary:

**Corollary 1** (Singular frameworks). *Let $d \in \mathcal{S}_0$. The Jacobian of the 2-cycles formation with restricted dynamics $F_r$ is generically of corank 1 for at least one framework attached to $d$.*

*Proof.* Recall that $A_e$ for the 2-cycles, given in Example 1 in Part I. A direct computation gives

$$J = u'(0) \begin{bmatrix} -z_1 \cdot z_1 & z_1 \cdot z_2 & 0 & 0 & -z_1 \cdot z_5 \\ 0 & -z_2 \cdot z_2 & z_2 \cdot z_3 & 0 & 0 \\ z_3 \cdot z_1 & 0 & -z_3 \cdot z_3 & 0 & z_3 \cdot z_5 \\ 0 & 0 & z_3 \cdot z_4 & -z_4 \cdot z_4 & 0 \\ -z_1 \cdot z_5 & 0 & 0 & z_4 \cdot z_5 & -z_5 \cdot z_5 \end{bmatrix}.$$

By Corollary 2, $u'(0) \neq 0$ generically. The first and last column are multiples of each other if $z_1$ is parallel to $z_5$, and it is easy to see that the first four columns are linearly independent. The corank is higher if, in addition, one of the $z_i$ is zero. ∎

We now prove Proposition 2.

*Proof of Proposition 2.* Consider a framework with $d \in \mathcal{S}_0$ and $z_1$ parallel to $z_5$. From Equation 27, we know that $\frac{\partial F_r}{\partial z}$ is generically of rank 4. Let $w$ be an eigenvector corresponding to the zero eigenvalue. We conclude using Lemma 9 that $w^T \frac{\partial F}{\partial d} = 0$. ∎

For the proof of Theorem 3 Bis, we exhibit a one-parameter family of systems by letting the value of $d$ vary around $\mathcal{S}_0$. An application of Sotomayor's theorem to this family of vector fields shows that around $\mathcal{S}_0$, a stable ancillary equilibrium is created.

*Proof of Theorem 3 Bis.* Recall that $\mathcal{F}_r = \{F(x) \mid u(0) = 0\}$. Fix $d_0 \in \mathcal{S}_0$ and write $d_0 = (d_1, d_2, d_3, d_4, d_5)$. We consider the *one parameter system* where only the distance between $x_1$ and $x_3$ is allowed to vary via change in the parameter $\mu$:

$$\begin{cases} \dot{x}_1 &= u(e_1)z_1 + u(e_5)z_5 \\ \dot{x}_2 &= u(e_2)z_2 \\ \dot{x}_3 &= u(z_3 \cdot z_3 - (d_3 + \mu))z_3 \\ \dot{x}_4 &= u(e_4)z_4 \end{cases} \quad (28)$$

with $\mu$ a scalar parameter. We will prove that conditions $(1), (2)$ and $(3)$ of Theorem 1 are generic for $\mathcal{F}$. Observe that we cannot use the second part of Sotomayor's result to deduce genericity from the three above conditions, since the set of vector fields we are looking at is much smaller than the one considered in Theorem 1.

From Corollary 1 and the fact that $u' \neq 0$ generically at a zero of $u$ (see Corollary 2), we know that the Jacobian of the 2-cycles at $\mathcal{S}_0$ has generically a unique zero eigenvalue for $F \in \mathcal{F}$. Hence condition $(1)$ is verified. Condition $(2)$ follows from Proposition 2.

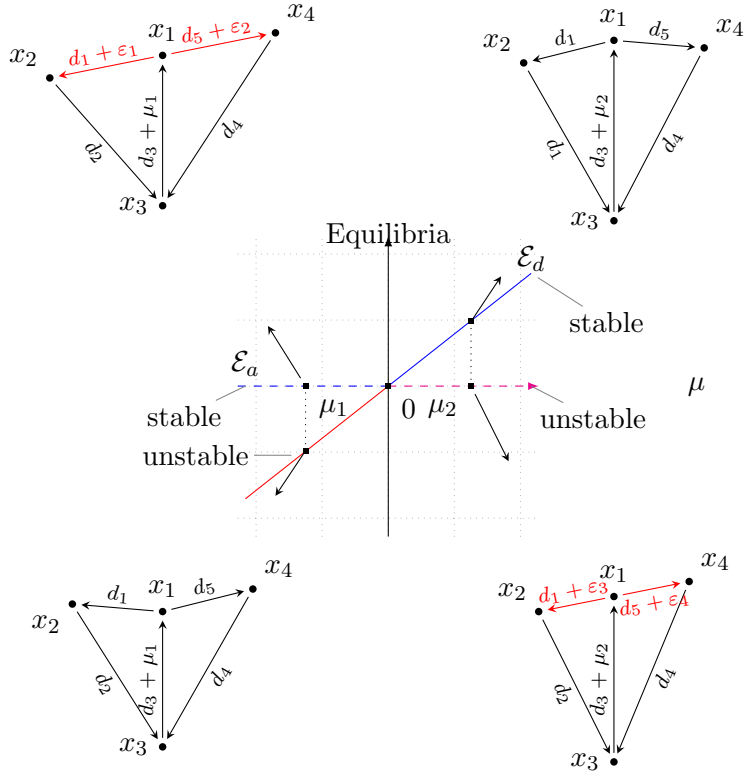

Figure 14: We illustrate the stability properties of ancillary and design equilibria around $\mathcal{S}_0$. Let the vector $(d_1, d_2, d_3, d_4, d_5) \in \mathcal{S}_0$. The horizontal dashed line corresponds ancillary equilibria and the slanted line that intersects it to design equilibria. They coincide at $\mu = 0$, as seen in Proposition 7.2; for $\mu \neq 0$ configurations in $\mathcal{S}_0$ are ancillary equilibria. For $\mu_1 < 0$, there is an ancillary equilibrium with $e_2, e_3, e_4 = 0$ but $e_1 = \varepsilon_1$ and $e_5 = \varepsilon_2$ and $z_1$ and $z_5$ aligned. It is illustrated in the top-left corner of the figure. This equilibrium is moreover stable. For $\mu_2 > 0$, there is a similar ancillary equilibria with $e_1 = \varepsilon_3$ and $e_5 = \varepsilon_4$, illustrated in the bottom-right corner, but this equilibrium is unstable. We see that around the bifurcation value $\mu_0$, there is a *transfer of stability* from $\mathcal{E}_d$ to $\mathcal{E}_a$. The orientation may be reversed (i.e. $\mu_1 > 0, \mu_2 < 0$ and all else the same in the figure) depending on the sign of the second derivatives in Theorem 3.

The second derivatives of condition (3) in Theorem 1 are, at a design equilibrium, the sum of two terms:

$$\frac{\partial^2 F}{\partial z^2} = u' \sum_l w_l v^T Q_l v + u'' \sum_l w_l v^T \tilde{Q}_l v$$

where $Q_l, \tilde{Q}_l \in \mathbb{R}^{10 \times 10}$ are obtained by evaluating the Hessian of $F_l$. One can easily check that $\sum_l w_l v^T Q_l v$ and $\sum_l w_l v^T \tilde{Q}_l v$ are generically non zero on $\mathcal{S}_0$. Set $a = \sum_l w_l v^T Q_l v$ and $b = \sum_l w_l v^T \tilde{Q}_l v$. We thus have $\frac{\partial^2 F}{\partial z^2}$ is zero only if

$$u'a + u''b = 0$$

when $u$ vanishes. Let $C \subset J^2(\mathbb{R}, \mathbb{R})$ be defined by the equations $au' + bu'' = 0$ and $u = 0$. Since $C$ is of codimension 2 in $J^2$, the 2-jet extension of $u$ is transversal to $C$ if and only if $u \neq 0$ or $au' + bu'' \neq 0$. We conclude using Thom's transversality theorem that $\frac{\partial^2 F}{\partial z^2}$ is generically non-zero.

Observe that $\frac{\partial F^2}{\partial z \partial d}$ is a constant matrix that does not depend on the configuration. Using a similar reasoning as above, we can conclude that $w^T \frac{\partial^2 F}{\partial z \partial d} v$ is generically non-zero. ∎

The proof of Theorem 4 does not assume the use of $\mathcal{F}_r$ over $\mathcal{F}$.

*Proof of Theorem 4.* We will show that there is a set of positive measure in $\mathcal{L}$ which cannot be made robustly type-A stable. We do so by showing that for any framework attached to distances in that set, the existence of a stable ancillary equilibrium is generic for $\mathcal{F}$

Denote by $\mathcal{S}^\varepsilon$ a tubular neighborhood of $\mathcal{S}$:

$$\mathcal{S}^\varepsilon = \{d \in \mathcal{L} \text{ s.t. } \exists d_0 \in \mathcal{S}, \text{ with } \|d - d_0\| < \varepsilon\}$$

and $\mathcal{S}_0^\varepsilon = \mathcal{S}^\varepsilon \cap \mathcal{S}_0$. The set $\mathcal{S}^\varepsilon$ contains frameworks where $z_1$ and $z_5$ are close to parallel. These frameworks are infinitesimally rigid and non-singular. Let $d \in \mathcal{S}_0^\varepsilon$ and $d_0 \in \mathcal{S}_l$ be such that there is $-\varepsilon < r < \varepsilon$ with $d = d_0 + (0, 0, \mu, 0, 0)$. Such $d_0$ and $\mu$ exist by definition of $\mathcal{S}_0^\varepsilon$.

Assume without loss of generality that $u$ is such that the design equilibria for $\mu > 0$ are stable, for the frameworks with either $x_1$ and $x_4$ on the same side or on opposite sides of $z_3$. If this is not the case, $\mathcal{E}_d \nsubseteq \mathcal{E}_s$ and the system is not type-A stable for that choice of $u$. Because for a generic $u$, the system undergoes a transcritical bifurcation when $\mu = 0$ by Theorem 3, and because $u' \neq 0$ generically when $u$ vanishes, we have that for $\varepsilon$ small enough, $\mathcal{E}_a$ contains the framework where $z_1$ is parallel to $z_5$ for all frameworks with $-\varepsilon < r < \varepsilon$. Furthermore, for either $\mu > 0$ or $\mu < 0$, we have that this framework is asymptotically stable, i.e. $\mathcal{E}_s \cap \mathcal{E}_a \neq \emptyset$. Hence, there is a set of positive measure of target frameworks in $\mathcal{S}_0^\varepsilon$ which contains a stable ancillary equilibrium and thus the system is not robustly type-A stable. ∎

# 8 Summary and conclusion

## 8.1 Summary

We now informally summarize the contents of this paper. The notions of type-A and strong type-A stability are extensions of the usual notion of global stability. These notions are akin to global stability in the sense that if a system is type A or strongly type-A stable, it will stabilize around a desired configuration almost surely; they extend global stability in that they acknowledge that the topology of the state-space of the system may introduce undesired equilibrium configurations. Under the assumption that the control system has a finite set $\mathcal{E}$ of equilibria, we expressed $\mathcal{E}$ as the disjoint union $\mathcal{E}_d \cup \mathcal{E}_a$ where $\mathcal{E}_d$ is the set of design equilibria—or the ones we would like the system to stabilize to—and $\mathcal{E}_a$ is the set of ancillary equilibria—or the equilibria that are unavoidable because of the topology of the space. A system is strongly type-A stable if all of the elements of $\mathcal{E}_d$ are locally stable, and no other equilibrium is stable; a system is type-A stable if only elements of $\mathcal{E}_d$, but not necessarily all of them, are stable. Both situations are of interest to formation control problems.

We established an algebraic framework to approach decentralized control problems and used it to obtain an obstruction to strong type-A stability for the 2-cycles. The framework showed in what sense one can claim that the individual agents do not have enough information to make the system strongly type-A stable. The development of this framework is the subject of current work.

We then proved that the 2-cycles cannot be made type-A stable. Our approach was to show that decentralized control laws of the type of Equation (14) exhibit a transcritical bifurcation generically, and from this concluded that ancillary stable equilibria exist generically.

A major technical aspect of this work is the one of genericity or, its counterpart, robustness. Loosely speaking, we say that a control law is *robust* if qualitative properties of $\mathcal{E}_a$ and $\mathcal{E}_d$ are not changed under small perturbations in the control. Consequently, the question may arise whether there are *non-robust* control laws, as impractical as they would be, that would render the system, say, type-A stable. These control laws would be such that either $u'(0) = 0$ or such that the second derivative of $F$ is zero, both situations being non-robust according to Theorem 3. In general, we believe that the answer to the existence of non-robust control laws that would make the system type-A stable is a negative one, because of the following ansatz: the more you increase the level of degeneracy at a bifurcation— e.g. by increasing its codimension [GH83, AAIS94]— the more complex the behavior of the system at the bifurcation is. Said differently, if one increases the level of degeneracy, the diagram of Figure 14 will increase in complexity in the sense that there will be a larger number of branches of equilibria compared to the two that are present in the current diagram. Hence, non-robust situations are more likely to exhibit a far more complex behavior and less likely to satisfy $\mathcal{E}_s \subset \mathcal{E}_d$.

The above argument also highlights the fact that, even though we argued for the "least amount of degeneracy", which implied robustness, we still have that the Jacobian of the 2-cycles is of corank 1 generically on $\mathcal{S}_0$. This unremoveable singularity is a *consequence of the distributed nature* of the control law, and thus is robust to perturbations of the control law. If one allows for more communication links, this situation may become non-generic and the transcritical bifurcation may then disappear.

Theorem 4 asserts the existence of a set of strictly positive measure in $\mathcal{L}$ such that the system is not robustly type-A stable. We can characterize this set more precisely. From the proof of Theorem 4, we see that for some frameworks around $\mathcal{S}$ the set of ancillary equilibria contains at least one stable element. These frameworks are, as is shown in Figure 14, on one side of $\mathcal{S}$ (since $\mathcal{S}$ is of codimension one, its sides are well-defined). A simple argument can be used to show that this region extends to frameworks quite distant from $\mathcal{S}$. Indeed, in order to observe a change in the type of equilibria for $\mu < 0$, one would need to go through another bifurcation. But from the form of system (14), one can see that no such bifurcations can happen until $\mu$ is such that we encounter another degeneracy of the framework (e.g. two or more agents aligned). One can show that target frameworks in the set

$$\mathcal{L}_c = \{d \in \mathcal{L}_0 \text{ s.t. one framework attached to } d \text{ has } x_1 \text{ in the convex hull of } x_2, x_3, x_4\}$$

result in systems which are not type-A table.

## 8.2 Conclusion

The main contributions of this paper can be summarized as follows:

- We introduced the notions of type-A and strong type-A stability.
- We introduced an algebraic framework to study decentralized control problems.
- We showed how the theory of singularities and bifurcations can be used in the context of type-A stability analysis.
- We applied all the concepts above to a particular decentralized system, namely the 2-cycles, to show that it was neither type-A nor strongly type-A stable.

The proofs of the results concerning type-A stability of the 2-cycles follow a common approach: we consider the set $\mathcal{F}$ of all possible feedback systems that respect the information flow and the invariance of formation control, as defined in part I, and show that no robust elements in this set satisfy necessary conditions for type-A stabilization.

The fact that a formation as simple as the 2-cycle cannot be stabilized globally, which may appear surprising to some, underscores the need to separate the constraint graph, which is to be handled by rigidity theory, from the information flow graph, for which the results of this paper introduce new approaches. In this context, we believe an answer to the two following questions would be of great interest to formation control:

- what is the minimal information flow graph of the 2-cycles that makes the system type-A stabilizable via a smooth feedback? Minimal here is understood as the least number of directed edges.

- does there exist non-continuous feedback that would stabilize the 2-cycles with the information flow graph studied in this paper?

We are currently developing the framework presented in this paper to provide general approaches and further our understanding of decentralized systems in the smooth and analytic cases.

## 9 Acknowledgments

We would like to thank Prof. B.D.O Anderson, Prof. R.W. Brockett, Prof. S. Morse as well as Alan O'Connor for helpful discussions. We are particularly grateful to Prof. Morse for introducing us to this problem.

## A  Example: $u(x) = x$

We illustrate the failure of the 2-cycles of Equation (14) to be type-A stable. We take $u(x) = x$ which is a straightforward extension of the control laws used in many works on directed formation control [CMY$^+$07, YADF09, KBF08].

The distributed system is then explicitly:

$$\begin{cases} \dot{x}_1 &= e_1 z_1 + e_5 z_5 \\ \dot{x}_2 &= e_2 z_2 \\ \dot{x}_3 &= e_3 z_3 \\ \dot{x}_4 &= e_4 z_4 \end{cases} \tag{29}$$

We show in Figure 15 the results of simulations for the vector of target distances $d_1 = 2.0, d_2 = 2.6, d_3 = 2.0, d_4 = 1.4, d_5 = 3.3$. We have that $d \in \mathcal{L}_c$, and we illustrate the formation such that $x_1$ is in the convex hull of $x_2, x_3$ and $x_4$ in Figure 15b. This configuration is unstable, whereas the configurations shown in Figures 15a and 15c are stable.

# B  Singularities of vector fields, jet spaces and transversality

The main tool handling genericity and robustness in function spaces is Thom's transversality theorem. We will arrive at the result by building onto the simpler concept of transversality of linear subspaces.

Thom's theorem roughly answers the following type of questions: given a function $u$ from a manifold $M$ to a manifold $N$, and some relations between the derivatives of different orders of this function (e.g. $u'' + u' - u = 0$), under what circumstances are these relations preserved under small perturbations of the function? For example, if a real-valued function has a zero at some point, under a small perturbation of this function, the zero will persist *generically* for $u$. On the other hand, if a real-valued function vanishes with its second derivative also being zero, under a small perturbation this property will be lost, see Figure 16. The crux of Thom's theorem is to show that considering only a "small subset" of perturbations (the integrable perturbations as we will see below) in the set of all perturbations in jet-spaces is sufficient.

Let $A, B \subset \mathbb{R}^n$ be linear subspaces. They are *transversal* if

$$\mathbb{R}^n = A \oplus B,$$

where $\oplus$ denotes the direct sum. For example, a plane and a line not contained in the plane are transversal in $\mathbb{R}^3$. The notion of transversality can be extended to maps as follows: given

$$F_1 : \mathbb{R}^n \to \mathbb{R}^m \text{ and } F_2 : \mathbb{R}^l \to \mathbb{R}^m,$$

we say that $F_1$ and $F_2$ are transversal at a point $(x_1, x_2) \in \mathbb{R}^n \times \mathbb{R}^l$ if one of the two following conditions is met:

1. $F_1(x_1) \neq F_2(x_2)$

2. If $F_1(x_1) = F_2(x_2)$, then the matrix $\begin{bmatrix} \frac{\partial F_1}{\partial x} \\ \frac{\partial F_2}{\partial x} \end{bmatrix}$ is of full rank.

In particular, if $l + n < m$ then $F_1$ and $F_2$ are transversal only where they do not map to the same point. This definition extends immediately to smooth functions between smooth manifolds: given

$$F_1 : M_1 \to N \text{ and } F_2 : M_2 \to N,$$

we say that $F_1$ and $F_2$ are transversal at $(x_1, x_2) \in M_1 \times M_2$ if either $F_1(x_1) \neq F_2(x_2)$ or $F(x_1) = F(x_2)$ and the tangent space of $N$ at $F(x_1)$ is the direct sum of the images of the tangent spaces of $M_1$ and $M_2$ under $F_1$ and $F_2$ respectively, i.e.

$$T_{F_1(x_1)} N = F_{1*} T_{x_1} M_1 \oplus F_{2*} T_{x_2} M_2$$

where $F_*$ is the push-forward [War83] of $F$.

**Example 6.** *Take $M_1 = \mathbb{R}$ and $M_2 = N = \mathbb{R}^2$ and let $F_1(x_1) = x_1 v + b_1$ and $F_2(x_2) = Ax_2 + b_2$, where $A \in \mathbb{R}^{2 \times 2}, b_2, v, x_2 \in \mathbb{R}^2$ and $x \in \mathbb{R}$. If $b_1 \neq b_2$, then $F_1(0) \neq F_2(0)$ and $F_1$ is transversal to $F_2$ at $(0,0,0)$. If $b_1 = b_2$, then $F_1(0) = F_2(0)$ and the functions are transversal if the span of $v$ and the columns of $A$ is $\mathbb{R}^2$.*

The notion of transversality that is of interest to us is a straightforward extension of the transversality of maps:

**Definition 6** (Transversality)**.** *Let $F : M \to N$ be a smooth map and let $C$ be a submanifold of $N$. Then $F$ is* transversal *to $C$ at a given point if, at that point, $F$ is transversal to the embedding $i : C \to N$ of $C$ into $N$.*

**Example 7.** *Take $N = \mathbb{R}^3$ with coordinates $u, v, w$ and $C$ be the u-v plane. Let $F : \mathbb{R} \to \mathbb{R}^3 : x \to [x, 2x, 3x]^T$. Then the map $F$ is transversal to $C$ everywhere.*

**Example 8.** *Let $F(x) : \mathbb{R} \to \mathbb{R}^3$ be any smooth curve in $\mathbb{R}^3$ and $C$ be the u-axis. At points where $F(x) \in C$, the tangent vector to $C$ and the tangent vector to $F$ will span at most a two-dimensional subspace in $\mathbb{R}^3$. Hence, $F$ is transversal to $C$ only at the points where $F(x) \notin C$.*

## B.1 Jet Spaces

Let $F, G : M \to N$ be smooth maps between smooth manifolds $M$ and $N$ endowed with a metric. We say that $F$ and $G$ are k-equivalent at $x_0 \in M$ if in a neighborhood of $x_0$ we have
$$\|F(x) - G(x)\| = o(\|x - x_0\|^k).$$
One can verify [Arn72] that k-equivalence is independent of the choice of metrics on $M$ and $N$ and that it is an equivalence relation on maps. In fact, the above definition can be recast as saying that $F$ and $G$ are 0-equivalent at $x_0$ if
$$F(x_0) = G(x_0),$$
1-equivalent if in addition
$$\frac{\partial F}{\partial x}\Big|_{x_0} = \frac{\partial G}{\partial x}\Big|_{x_0},$$
and so forth. We define the k-jet of a smooth map to be its k-equivalence class:

**Definition 7.** *The k-jet of $F : M \to N$ at $x_0$ is*
$$J_{x_0}^k(F) = \{G : M \to N \text{ s.t. } G \text{ is k-equivalent to } F\}.$$

Hence, the 0-jet of $F$ at $x_0$ is $F(x_0)$; the 1-jet is $(F(x_0), \frac{\partial F}{\partial x}|_{x_0})$, etc. For example, the constant function 0 and $\sin(x)$ have the same 0-jet at $x = 0$ and $x$ and $\sin(x)$ have the same 1-jet at 0.

We define:
$$J^k(M, N) = \text{ Space of k- jets from } M \text{ to } N.$$

A 0-jet is thus determined by a point in $M$ and a point in $N$, and thus $J^0(M, N)$ nothing more than the Cartesian product of $M$ and $N$:
$$J^0(M, N) = M \times N.$$

Since a 1-jet is determined by a pair of points, for the 0-jet part, and a matrix of dimension $\dim M \times \dim N$, for the Jacobian of the function at $x_0$, we see that $\dim J^1(M, N) = \dim M + \dim N + \dim M \dim N$. We cannot say in general that $J^1(M, N)$ is the cartesian product of $J^0$ with $\mathbb{R}^{m \times n}$ since the product may be twisted. Similar relations are obtained for higher jet-spaces [Arn72]

Given a function $F : M \to N$, we call its *k-jet extension* the map given by
$$j_F^k(x) : M \to J^K(M, N) : x \to (F(x), \frac{\partial F}{\partial x}(x), \ldots, \frac{\partial^k F}{\partial x^k}(x)).$$

**Example 9.** *Let $M = N = \mathbb{R}$. The k-jet space is $J^k(\mathbb{R}, \mathbb{R}) = \mathbb{R} \times \mathbb{R} \times \ldots \times \mathbb{R} = \mathbb{R}^{k+2}$. Take $F(x) = \sin(x)$; the 2−jet extension of $F$ is*
$$j_{\sin}^2(x) = (x, \sin(x), \cos(x), -\sin(x)).$$
*If we take $M = N = \mathbb{R}^2$ and $F(x) = Ax$ for $A \in \mathbb{R}^{2 \times 2}$, then $j_{Ax}^k(x) = (x, Ax, A, 0, \ldots, 0)$.*

**Remark 2.** *The concepts presented in this section also trivially apply to vector fields on $M$, by letting $N = TM$.*

While to any function $F : M \to N$, we can assign a k-jet extension $j_F^k : M \to J^k(M, N)$, the inverse is not true: there are maps $G : M \to J^k(M, N)$ which do not correspond to functions from $M$ to $N$ as there are some obvious integrability conditions that need to be satisfied. For example, if we let
$$G : \mathbb{R}^n \to J^1(\mathbb{R}^n) : G(x) = (x, Ax, B),$$
then $G$ is a 1-jet extension of a function if and only if $B = A$.

The power of the transversality theorem of Thom is that it allows one to draw conclusions about transversality properties in general, and genericity in particular, by *solely looking at perturbations in jet spaces that are jet extensions*—a much smaller set than all perturbations in jet-spaces, since these include the much larger set of non-integrable perturbations.

We recall that the $C^r$ topology is a metric topology. It is induced by a metric that takes into account the function and its first $r$ derivatives. We have:

**Theorem 5** (Thom's transversality). *Let $C$ be a regular submanifold of the jet space $J^k(M, N)$. Then the set of maps $f : M \to N$ whose k-jet extensions are transversal to $C$ is an everywhere dense intersection of open sets in the space of smooth maps for the $C^r$ topology, $1 \leq r \leq \infty$.*

A typical application of the theorem is to prove that vector fields with degenerate zeros are not generic. We here prove a version of this result that is of interest to us.

**Corollary 2.** *Functions in $\mathcal{C}^\infty(M)$ whose derivative at a zero vanish are not generic.*

In other words, the corollary deals with the intuitive fact that if $u(x) = 0$, then generically $u'(x) \neq 0$. This result also goes under the name of weak-transversality theorem [GG74].

*Proof.* Consider the space of 0-jets $J^0(M, \mathbb{R})$. In this space, let $C$ be the set of 0-jets which vanish, i.e. $C = (x, 0) \subset J^0$. A function $u$ is transversal to this set if either it does not vanish, or where it vanishes we have that the matrix

$$\begin{bmatrix} 1 & 1 \\ 0 & \frac{\partial f}{\partial x} \end{bmatrix}$$

is of full rank. Hence, transversality to $C$ at a zero implies that the derivative of the function is non-zero. The result is thus a consequence of Theorem 5. ∎

To picture the situation geometrically, recall that $J^0(M, \mathbb{R})$ is simply $M \times \mathbb{R}$. Hence $C$ is $M \times 0 \subset J^0(M, \mathbb{R})$. The result says that any function that intersects $C$ without crossing (and hence with a zero derivative) will, under a generic perturbation, either cross $C$ or not intersect $C$ at all, since these two eventualities result in transversality. Figure 16 provides an illustration when $M = \mathbb{R}$.

## C Factorization Lemma

LEMMA 3. (Factorization Lemma) *Let $\bar{x} = [\bar{x}_2, \bar{x}_3, \bar{x}_4]$ represent an equilibrium framework for the 2-cycles. For the dynamics of Equation (14), $F(\bar{x})$ is given by*

$$-F(\bar{x}) = \det(J) = p(\bar{x})q(\bar{x})$$

*where*

$$q(\bar{x}) = (k_2 k_3 k_4)(k_{11} k_{52} - k_{12} k_{51})$$

*and*

$$p(\bar{x}) = -\det(A_1)\det(A_2)\det(A_3)\det(A_4)$$

where the $A_i$ are $2 \times 2$ matrices given by

$$A_1 = \begin{bmatrix} | & | \\ \bar{z}_1 & \bar{z}_3 \\ | & | \end{bmatrix}, A_2 = \begin{bmatrix} | & | \\ \bar{z}_1 & \bar{z}_5 \\ | & | \end{bmatrix}, A_3 = \begin{bmatrix} | & | \\ \bar{z}_3 & \bar{z}_4 \\ | & | \end{bmatrix}, A_4 = \begin{bmatrix} \bar{z}_1 \cdot \bar{z}_3^{\perp} & \bar{z}_4 \cdot \bar{z}_3^{\perp} \\ \bar{z}_2 \cdot \bar{z}_2 & \bar{z}_4 \cdot \bar{z}_4 \end{bmatrix}.$$

*Proof.* We have that when $d_3 = 0$, $\det(J) = \det(A_1) = 0$ and the Lemma is verified. We thus assume that $d_3 \neq 0$ and rotate the framework so that $\bar{x}_3 = [0, -d_3]$ as in Figure 7b. Using these coordinates, we show that

$$\det(J) = (k_2 k_3 k_4)(k_{11} k_{52} - k_{12} k_{51}) p_1 p_2 p_3 p_4 \tag{30}$$

where

$$\begin{cases} p_1(\bar{x}_2, \bar{x}_4) &= d_3 \bar{x}_{21} \\ p_2(\bar{x}_2, \bar{x}_4) &= d_3 \bar{x}_{41} \\ p_3(\bar{x}_2, \bar{x}_4) &= (\bar{x}_{21} \bar{x}_{42} - \bar{x}_{22} \bar{x}_{41}) \\ p_4(\bar{x}_2, \bar{x}_4) &= d_3 \bar{x}_{41} d_2^2 - d_3 \bar{x}_{21} d_4^2 \end{cases} \tag{31}$$

We have that $p_i = \det(A_i)$ when $A_i$ is expressed in the coordinates described above. As a consequence of the invariance of $\det(J)$ under rotations, it is thus enough to prove that Equations (30) and (31) hold.

We express $J$, given in Corollary 1, in these coordinates:

$$J = \begin{bmatrix} -\bar{x}_{21}(k_{11}\bar{x}_{21} + k_{51}\bar{x}_{41}) - \bar{x}_{22}(k_{11}\bar{x}_{22} + k_{51}\bar{x}_{42}) & -k_2 x_{21}^2 - k_2 \bar{x}_{22}(x_{22} + d_3) & 0 \\ 0 & -k_2 x_{21}^2 - k_2(x_{22} + d_3)^2 & -d_3 k_3 (x_{22} + d_3) \\ d_3(k_{11}\bar{x}_{22} + k_{51}\bar{x}_{42}) & 0 & -d_3^2 k_3 \\ 0 & 0 & -d_3 k_3 (x_{42} + d_3) \\ -\bar{x}_{41}(k_{11}\bar{x}_{21} + k_{51}\bar{x}_{41}) - \bar{x}_{42}(k_{11}\bar{x}_{22} + k_{51}\bar{x}_{42}) & 0 & 0 \end{bmatrix}$$

$$\begin{bmatrix} 0 & -\bar{x}_{21}(k_{52}\bar{x}_{41} + k_{12}\bar{x}_{21}) - \bar{x}_{22}(k_{52}\bar{x}_{42} + k_{12}\bar{x}_{22}) \\ 0 & 0 \\ 0 & d_3(k_{52}\bar{x}_{42} + k_{12}\bar{x}_{22}) \\ -k_4 x_{41}^2 - k_4(x_{42} + d_3)^2 & 0 \\ -k_4 x_{41}^2 - k_4 \bar{x}_{42}(x_{42} + d_3) & -\bar{x}_{41}(k_{52}\bar{x}_{41} + k_{12}\bar{x}_{21}) - \bar{x}_{42}(k_{52}\bar{x}_{42} + k_{12}\bar{x}_{22}) \end{bmatrix}$$

where we recall that

$$k_i(d_i) = \frac{\partial}{\partial x} u_i(d_i; x)|_{x=0}$$

$i = 2, 3, 4$ and

$$k_{i1} = \frac{\partial}{\partial x} u_i(d_1, d_5; x, y, z_1 \cdot z_5)|_{x,y=0}, \quad k_{i2} = \frac{\partial}{\partial y} u_i(d_1, d_5; x, y, z_1 \cdot z_5)|_{x,y=0}$$

for $i = 1, 5$.

We first factor $k_i$, $i = 2, 3, 4$ from the second, third and fourth column respectively and $d_3$ from the third column. We use the notation $J_{IL}$ where $I, L$ are subsets of $\{1, 2, \ldots, n\}$ to refer to the submatrix of $J$ with rows and columns indexed in $I$ and $L$ respectively. We evaluate the determinant by expanding it along its second column:

$$\det(J) = -d_3 k_2 k_3 k_4 \left[ (\bar{x}_{21}^2 + \bar{x}_{22}(\bar{x}_{22} + d_3)) \det(J_{2345,1345}) \right.$$
$$\left. - (\bar{x}_{21}^2 + (\bar{x}_{22} + d_3)^2) \det(J_{1345,1345}) \right]$$

After some algebraic manipulations, we obtain

$$\det(J_{2345,1345}) = d_3 \bar{x}_{41}(k_{11}k_{52} - k_{12}k_{51})(\bar{x}_{21}\bar{x}_{42} - \bar{x}_{22}\bar{x}_{41})((\bar{x}_{42} + d_3)^2 + \bar{x}_{41}^2)(\bar{x}_{22} + d_3)$$
$$\det(J_{1345,1345}) = -d_3 \bar{x}_{41}(k_{11}k_{52} - k_{12}k_{51})(\bar{x}_{21}\bar{x}_{42} - \bar{x}_{22}\bar{x}_{41})(\bar{x}_{22}((\bar{x}_{42} + d_3)^2 + \bar{x}_{41}^2) + d_3 \bar{x}_{21}\bar{x}_{41}).$$

Pulling out the common factors, we have

$$\det(J) = -d_3^2 k_1 k_2 k_3 \bar{x}_{41} \bar{x}_{21}(k_{11}k_{52} - k_{12}k_{51})(\bar{x}_{21}\bar{x}_{42} - \bar{x}_{22}\bar{x}_{41}) \left[ (\bar{x}_{21}^2 + \bar{x}_{22}(\bar{x}_{22} + d_3)) \right.$$
$$((\bar{x}_{42} + d_3)^2 + \bar{x}_{41}^2)(\bar{x}_{22} + d_3) - (\bar{x}_{21}^2 + (\bar{x}_{22} + d_3)^2) (\bar{x}_{22}((\bar{x}_{42} + d_3)^2 + \bar{x}_{41}^2) + d_3 \bar{x}_{21}\bar{x}_{41}) \Big]$$

In the term between brackets in the previous equation, the terms $\bar{x}_{21}^2(\bar{x}_{22}((\bar{x}_{42}+d_3)^2+\bar{x}_{41}^2))$ and $(\bar{x}_{22} + d_3)^2(\bar{x}_{22}((\bar{x}_{42} + d_3)^2 + \bar{x}_{41}^2)$ simplify and there is left

$$d_3 \bar{x}_{21}^2((\bar{x}_{42} + d_3)^2 + \bar{x}_{41}^2) - d_3 \bar{x}_{21}^3 \bar{x}_{41} - (\bar{x}_{22} + d_3)^2 d_3 \bar{x}_{21} \bar{x}_{41}$$
$$= d_3 \bar{x}_{21} \left( \bar{x}_{21}((\bar{x}_{42} + d_3)^2 + \bar{x}_{41}^2) - d_3 \bar{x}_{41}((\bar{x}_{22} + d_3)^2 + \bar{x}_{21}^2) \right) = -p_4$$

since $d_2^2 = (\bar{x}_{22} + d_3)^2 + \bar{x}_{21}^2$ and $d_4^2 = (\bar{x}_{42} + d_3)^2 + \bar{x}_{41}^2$. ∎

## D  Proof of Theorem 3

We extend the results that were proved in Section 7. The proof of Theorem 3 follows the lines of the proof of Theorem 3bis. We will prove that the conclusions of Proposition 2 and Corollary 1 hold in the case of the more general dynamics given by

$$F(z) = \begin{bmatrix} u_2(d_2; e_2)z_2 - u_1(d_1, d_5; e_1, e_5, z_1 \cdot z_5)z_1 - u_5(d_1, d_5; e_1, e_5, z_1 \cdot z_5)z_5 \\ u_3(d_3; e_3)z_3 - u_2(d_2; e_2)z_2 \\ u_1(d_1, d_5; e_1, e_5, z_1 \cdot z_5)z_1 + u_5(d_1, d_5; e_1, e_5, z_1 \cdot z_5)z_5 - u_3(d_3; e_3)z_3 \\ u_3(d_3; e, 2)z_3 - u_3(d_4; e_4)z_4 \\ u_4(d_4; e_4)z_4 - u_1(d_1, d_5; e_1, e_5, z_1 \cdot z_5)z_1 - u_5(d_1, d_5; e_1, e_5, z_1 \cdot z_5)z_5 \end{bmatrix} \quad (32)$$

and set

$$\mathcal{F} = \{F(z) \mid u_i \in \mathcal{U}\}.$$

PROPOSITION 2. *Let $d \in \mathcal{S}_0$. There is a non-zero vector $w \in \mathbb{R}^{10}$ such that $w^T \frac{\partial F}{\partial z}|_{e_i=0,d} = w^T \frac{\partial F}{\partial d}|_d = 0$ for at least one framework attached to d.*

*Proof.* We start by evaluating the differential $\frac{\partial F}{\partial d}$. We have

$$\frac{\partial F_1}{\partial d_1} = -(\frac{\partial u_1}{\partial d_1} + \frac{\partial u_1}{\partial x})z_1 - (\frac{\partial u_5}{\partial d_1} + \frac{\partial u_5}{\partial x})z_5$$

and similar relations for the other entries $\frac{\partial F_i}{\partial d_j}$. Recall the definition of $z_i'$ in Proposition 4 in Part I: $z_i' = (u_x z_i + u_y z_j)$ if $z_i$ originates from an agent with two coleaders given by $z_i$ and $z_j$, and $z_i' = u_x z_i$ if $z_i$ originates from an agent with a single coleader. We define

$$\begin{cases} z_1'' &= z_1' + \frac{\partial u_1}{\partial d_1}z_1 + \frac{\partial u_5}{\partial d_1}z_5 \\ z_2'' &= z_2' + \frac{\partial u_2}{\partial d_2}z_2 \\ z_3'' &= z_3' + \frac{\partial u_3}{\partial d_3}z_3 \\ z_4'' &= z_4' + \frac{\partial u_4}{\partial d_4}z_4 \\ z_5'' &= z_5' + \frac{\partial u_1}{\partial d_5}z_1 + \frac{\partial u_5}{\partial d_5}z_5 \end{cases} \quad (33)$$

Define $Z''$ as in Equation 15 with $z_i$ replaced by $z_i''$. Some simple algebraic manipulations yield

$$\frac{\partial F}{\partial d} = A_e Z''^T. \quad (34)$$

We claim that when $z_1$ is parallel to $z_5$, we can find diagonal matrices $D_1$ and $D_2$ such that

$$Z' = D_1 Z = D_2 Z''. \quad (35)$$

Indeed, for $i = 2, 3, 4$, it is immediate from the definitions of $z_i$, $z_i'$ and $z_i''$ that there exists $\alpha, \beta_i \in \mathbb{R}$ such that

$$z_i = \alpha_i z_i' = \beta_i z_i'', i = 2, 3, 4.$$

This equation does not hold in general for $i = 1, 5$ unless there exists $\gamma \in \mathbb{R}$ such that $z_1 = \gamma z_5$. The $\alpha_i, \beta_i$ are the entries of $D_1$ and $D_2$. Recall that according to Proposition 4 in Part I, we have

$$\frac{\partial F}{\partial z} = A_e Z'^T Z \quad (36)$$

Putting Equations (34), (35), (36) together, we obtain that

$$\frac{\partial F}{\partial d} D_2 Z = \frac{\partial F}{\partial Z}.$$

If all the $z_i$ are non-zero, then $D_2 Z$ is of full rank and we conclude that

$$w^T \frac{\partial F}{\partial d} = 0 \Leftrightarrow w^T \frac{\partial F}{\partial z} = 0. \quad (37)$$

∎

COROLLARY 1 [Singular frameworks]. *Let $d \in \mathcal{S}_0$. The Jacobian of the 2-cycles formation is generically of corank 1 for at least one framework attached to d.*

*Proof.* A direct computation using Corollary 1 in part I and the edge-adjacency matrix of the 2-cycles gives

$$J = ZA_e Z'^T = \begin{bmatrix} -z_1 \cdot z_1' & z_1 \cdot z_2' & 0 & 0 & -z_1 \cdot z_5' \\ 0 & -z_2 \cdot z_2' & z_2 \cdot z_3' & 0 & 0 \\ z_3 \cdot z_1' & 0 & -z_3 \cdot z_3' & 0 & z_3 \cdot z_5' \\ 0 & 0 & z_3 \cdot z_4 & -z_4 \cdot z_4 & 0 \\ -z_1 \cdot z_5' & 0 & 0 & z_4 \cdot z_5' & -z_5 \cdot z_5' \end{bmatrix}.$$

By Corollary 2, $u_i'(0) \neq 0$ generically, hence $z_i'$ are generically non-zero. For the framework attached to $d$ in $\mathcal{S}_0$ such that $z_1$ is parallel to $z_5$, the first and last column are multiples of each other, and it is easy to see that the first four columns are linearly independent. The corank is higher if, in addition, one of the $z_i$ is zero. ∎

The proof of Theorem 3 is similar to the one of Theorem 3bis from this point on.

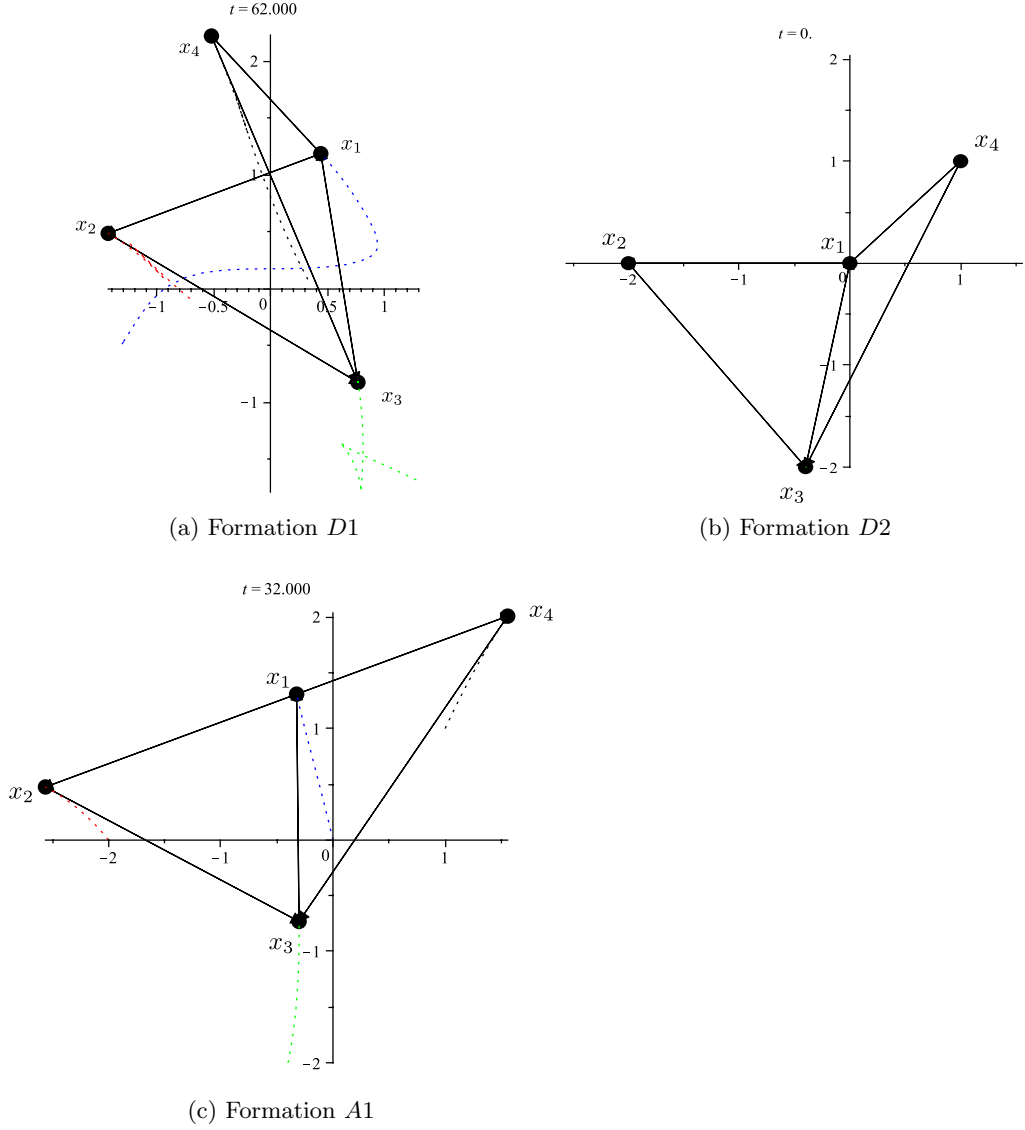

(a) Formation $D1$

(b) Formation $D2$

(c) Formation $A1$

Figure 15: Simulation results for the decentralized system of Equation (29) with a $d \in \mathcal{L}_c$. The dotted lines represent the trajectories followed by the agents. Two of the configurations in $\mathcal{E}_d$ are depicted in Figures 15a and 15b. The other two configurations in $\mathcal{E}_d$ are their mirror symmetric. An ancillary equilibrium is shown in Figure 15c. A linearization of the system gives that the spectra of the Jacobians are given by $(-17.5 + 1.3i, -17.5 - 1.3i, -11.9, -7.9, -0.6)$, $(0.6, -18.6 + 3i, -18.6 - 3i, -9.4 + 3.1i, -9.4 - 3.1i)$ and $(-23.4 + 4.8i, -23.4 - 4.8i, -11 + 2.8i, -11 - 2.8i, -1.6)$ for the formations $D1$, $D2$, and $A1$ respectively. Hence $A1$ is locally exponentially stable and the system is not type-A stable.

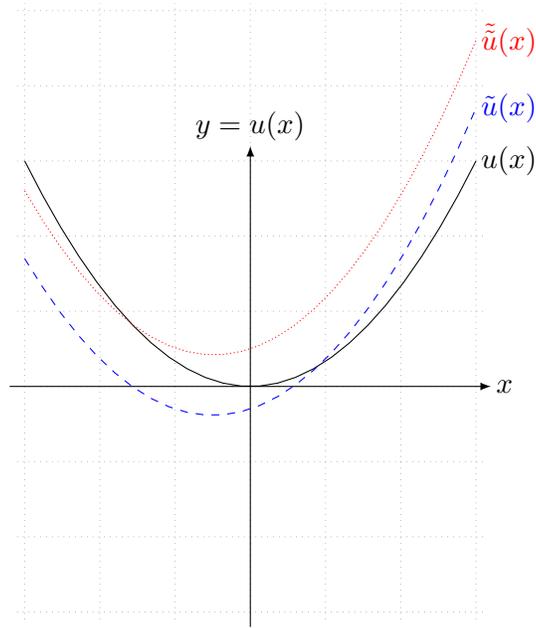

Figure 16: If we let $\mathcal{P}$ be the property of vanishing with a zero derivative. We will prove in this section that $\rceil\mathcal{P}$ is generic and thus $\mathcal{P}$ is not robust. Let $u(x)$ be a function which satisfy $\mathcal{P}$. For any small perturbations, it will either vanish with a non-zero derivative—as illustrated with $\tilde{u}(x)$, dashed curve— or not vanish at all—as illustrated with $\tilde{\tilde{u}}(x)$, dotted curve. Both $\tilde{u}(x)$ and $\tilde{\tilde{u}}(x)$ are transversal to the manifold defined by $y = 0$ everywhere, whereas $u(x)$ is not.